\documentclass{article}

\usepackage[preprint, nonatbib]{neurips_2023}

\usepackage[utf8]{inputenc} 
\usepackage[T1]{fontenc}    
\usepackage{hyperref}       
\usepackage{url}            
\usepackage{booktabs}       
\usepackage{amsfonts}       
\usepackage{nicefrac}       
\usepackage{microtype}      
\usepackage{xcolor}         
\usepackage{biblatex}
\usepackage{algpseudocode}
\usepackage{amsthm}
\usepackage{amsmath}
\usepackage{amssymb}
\usepackage{graphicx}
\usepackage{caption}
\usepackage{subcaption}
\usepackage{comment}
\usepackage{tikz-cd}
\usetikzlibrary{calc}
\usetikzlibrary{cd}

\bibliography{references.bib}

\newtheorem{theorem}{Theorem}[section]

\theoremstyle{definition}
\newtheorem{definition}[theorem]{Definition}
\theoremstyle{remark}
\newtheorem{remark}[theorem]{Remark}

\newcommand*\samethanks[1][\value{footnote}]{\footnotemark[#1]}
\newcommand{\Set}{\mathrm{Set}}
\newcommand{\Gr}{\mathrm{Gr}}

\newenvironment{rcases}
  {\left.\begin{aligned}}
  {\end{aligned}\right\rbrace}

\title{Applying language models to algebraic topology: generating simplicial cycles using multi-labeling in Wu's formula}

\author{
  \kern-1.5em Kirill~Brilliantov\thanks{Joint first author} \\
  \kern-1.5em Constructor University\\
  \And
  Fedor Pavutnitskiy\samethanks  \\
  BIMSA\thanks{Yanqi Lake Beijing Institute of Mathematical Sciences and Applications} \\
  \AND
  Dmitry Pasechnyuk \\
  MBZUAI, MIPT, IITP RAS\\
  \And
  German Magai \\
  Higher School of Economics \\
}

\begin{document}

\maketitle

\begin{abstract}

Computing homotopy groups of spheres has long been a fundamental objective in algebraic topology. Various theoretical and algorithmic approaches have been developed to tackle this problem. In this paper we take a step towards the goal of comprehending the group-theoretic structure of the generators of these homotopy groups by leveraging the power of machine learning. Specifically, in the simplicial group setting of Wu's formula, we reformulate the problem of generating simplicial cycles as a problem of sampling from the intersection of algorithmic datasets related to Dyck languages. We present and evaluate language modelling approaches that employ multi-label information for input sequences, along with the necessary group-theoretic toolkit and non-neural baselines.

\end{abstract}

\section{Introduction}\label{sec:Introduction}

Mathematical progress often depends on the identification of patterns and the formulation of conjectures, which are statements believed to be true but yet to be proven. For centuries, mathematicians have relied on data, from early prime tables to modern computer-generated data, to aid in this process. While computational techniques have proven useful in some areas of pure mathematics~\cite{appel19774color}, \cite{hales2005kepler}, the full potential of artificial intelligence (AI) in the field is yet to be fully explored.

In mathematics, AI has demonstrated remarkable abilities in tasks such as identifying counterexamples to conjectures~\cite{wagner2021constructions}, accelerating calculations~\cite{peifer2020learning}, generating symbolic solutions~\cite{lample2019deep}, detecting structural patterns in mathematical objects~\cite{he2021machinelearning}, automated generation of formal proofs~\cite{han2021proof}, \cite{polu2022formal} and guiding human intuition to making hypothesis~\cite{davies2021guiding}. These capabilities offer vast potential for expanding mathematical knowledge and advancing the field.

Inspired by these recent advanced in AI driven mathematical research, we propose a proof-of-concept framework for generating simplicial cycles in simplicial groups. This serves as an initial step towards a broader initiative of utilizing machine learning techniques to
sample elements from the homotopy groups of spheres and general spaces. In the current work\footnote{Code available at \url{https://github.com/ml-in-algebraic-topology/gen-simplicial-cycles}} we apply our ideas in particular case of free simplicial group which models (loops over) a two-dimensional sphere. We adopt a simplicial and group-theoretic approach
due to its combinatorial nature and computer-friendly representation.

The theoretical background for Wu's formula~\eqref{eq: Wu formula} is outlined in Section~\ref{sec: group theory} with a more detailed exposition given in Appendix~\ref{apdx: background on simplicial homotopy theory}. This formula expresses our object of interest as a certain quotient of the intersection of normal subgroups of free group. Generating elements from this intersection $\cap_i R_i$ (see notation~\eqref{eq: definition of R_i} below) is the main topic of the present work. The key features of the problem include the similarity between $R_i$'s and sets of Dyck-$n$ words (allowing a suitable distribution on it to sample from) and a clear algorithmic procedure for determining if a given element $w$ lies within the intersection.

We develop the following non-neural-based baselines for this problem: random search, evolutionary method and greedy algorithm. These baselines are extensively described in Section~\ref{subsec: baselines}. Concurrently, our main deep learning approaches described in Section~\ref{subsec: model description} are based on language modelling and differ in the various usages of training dataset's meta information.

The goal of introducing the deep learning approaches is to increase the performance in terms of the number of generated words from the full intersection, while also achieving scalability of the model. This is crucial due to the rapid growth in computational complexity observed in the baseline algorithms, rendering them impractical for application in higher dimensions. The relevant set of metrics, datasets description and overall evaluation results are presented in Section~\ref{sec: evaluation} together with some relevant discussions pertaining to the findings.

\section{Related work}\label{sec:background and related work}

The computation of homotopy groups of spheres stands as one of the most challenging problems in mathematics.  Over the past century, the complexity of this problem has spurred the development of numerous sophisticated mathematical techniques. For instance, in the 1950s, Toda introduced the concept of secondary operations (now called Toda brackets), which allows him to compute the homotopy group of spheres up to a certain dimension~\cite{toda1963compositional}. In the 1960s, unstable Adams spectral sequences were developed, providing a general framework for
such computations and extending the range of results achieved by Toda~\cite{adams1958spectralsequences}. In the 1990s a fresh approach using Goodwillie tower of identity functor~\cite{arone1999goodwillie}, \cite{behrens2012goodwillie} was employed, leveraging modern techniques to recompute Toda’s range.

In recent years, advancements have been made in computational and computer
science approaches concerning the computation of homotopy groups.

A classical result by Brown~\cite{brown1957finitelycomputability} established that homotopy groups of spheres are finitely computable. However, it should be noted that the algorithm presented in Brown’s work, as discussed in~\cite{cadek2014polynomialtime}, exhibits superexponential complexity. While this result lays a solid theoretical foundation for the development of algorithms and computational methods for computing these groups, the proof itself does not directly yield practical algorithms for their computation. Nonetheless, it has sparked new avenues of research and possibilities for researchers interested in devising effective computational methods for computing homotopy groups.

While the problem of computing even the weaker invariants like the rationalization $\pi_n X \otimes \mathbb{Q}$ for a general space $X$ is known to be $\#P$-hard~\cite{anick1987computational}, there have been notable advancements in algorithmic computation of homotopy groups and homology of spaces.

The authors of the Kenzo system~\cite{sergeraert1994computability},\cite{rubio2002constructive}, \cite{romero2006computing} since the 1990s have been actively developing an algorithmic framework based on spectral sequences to compute homotopy groups of certain classes of spaces. Their concept of \emph{spaces with effective homology} was further enhanced in~\cite{cadek2014polynomialtime}, where a polynomial-time algorithm was introduced for computing $\pi_n X$ for a \textbf{fixed} value of $n$ within the framework of \emph{simplicial sets with polynomial time homology}.

Besides the computational challenges, one limitation of current approaches is that $\pi_n X$ is typically computed as an abelian group, without providing a detailed account of its internal homotopical structure, such as a description of generators and relations in algebraic, geometric, or topological terms. However, there have been efforts~\cite{filakovsky2018computing} proposing approaches to address this limitation.

The utilization of machine learning techniques in the calculation of homotopy groups of spaces is a relatively novel and emerging approach. So far, the corresponding research concerning both machine learning and algebraic topology, primarily dedicated to enhancing machine learning algorithms with ideas from algebraic topology, but lacks exploration of the inverse. The present work initiates the process of filling this gap.

\section{Group theory of Wu's formula}\label{sec: group theory}
To keep the preliminary exposition concise, we will provide a brief overview of the algebraic constructions employed in this paper, reserving more formal and detailed information for Appendix~\ref{apdx: background on group theory}. Wu's formula for $\pi_{n}S^2$, introduced in~\cite{wu2001combinatorial}, combines the simplicity of group theoretic language of free group with a powerful combinatorial structure of simplicial groups. In our perspective, the group theoretic structures in Wu's formula, related to formal languages, presents a compelling target for the application of language models. Additionally, the usage of simplicial groups opens up further avenues for incorporating machine learning to a more sophisticated mathematical constructions of simplicial homotopy theory, like the unstable Adams spectral sequence~\cite{bousfield1966mod}.

By \emph{free group} $F$ on a set ${x_1,\dots x_n}$ we will mean a set of words in letters $x_i$ (called \emph{generators}) and their inverses with group operation of string concatenation. The identity $x_ix_i^{-1} = x_i^{-1}x_i = 1$ is assumed in $F$,
here $1$ is an empty word. Additionally, we will denote by $x_0 = x_1\dots x_n$ the product of all generators $x_i$. We will further use the following notation: $u^v = v^{-1}uv$ called \emph{conjugation} and $[u,v] = u^{-1}v^{-1}uv$ called \emph{commutator} for $u,v \in F$.

We will be interested in generating words from the following subsets of $F$:
\begin{equation}\label{eq: definition of R_i}
  R_i = \{x_i^{\pm y_1}\dots x_i^{\pm y_k} \, | \,k\in\mathbb{N},\, y_i \in F\}, \ i = 0..n.  
\end{equation}

The subgroups $R_i$, also denoted by $\langle x_i \rangle^F$, to emphasize the generating element, are called the \emph{normal closures} of $x_i$'s. These subgroups can be conceptually understood as follows: for $i > 0$ the subgroup $R_i$ consists of words that become trivial after substituting $x_i$ for $1$, for example:
\[
[x_1, x_2]x_1 \in R_1, \ [x_1, x_2]x_1 |_{x_1 = 1} = 1.
\]
Similarly, $R_0$ consist of words which become trivial after any of the equivalent substitutions like $x_1\mapsto x_n^{-1}x_{n-1}^{-1}\dots x_2^{-1}$, which maps $x_0$ to $1$. Equivalently, $R_0$ consists of words which become trivial upon removal the cyclic permutations of $x_0$.

As a result, the intersection $\cap_{i = 0}^n R_i$ consists precisely of words that become trivial when any of the (cyclic permutations of) $x_i$ is substituted with $1$. The objective of this paper is to propose a several deep learning generators of such words.

Since $[u,1] = [1,v] = 1$ for all $u,v \in F$, it is clear that the elements of the form $[x_i,x_j]$ are in the intersection $R_i \cap R_j$ and, more generally, $[R_i,R_j] \subseteq R_i \cap R_j$. This idea can be iterated by introducing the \emph{iterated commutators} 
\begin{equation}\label{eq: iterated commutator notation}
[u_1,\dots, u_k] = [[[u_1,u_2],u_3]\dots,u_k], \ u_i \in F,
\end{equation}

and for a set of subgroups $\{R_{i_1},\dots, R_{i_k}\}$ we have: 
\begin{equation}\label{eq: inclusion of symm comm into intersection}
    \prod_{\sigma \in \Sigma_k} [R_{\sigma(i_1)},\dots R_{\sigma(i_k)}] \subseteq \bigcap_{j=1}^k R_{i_j},
\end{equation}
where the left-hand side have a product over all permutations on $k$ letters, reflecting the symmetric nature of the intersection on the right-hand side. We will abbreviate the left-hand side of the inclusion~\eqref{eq: inclusion of symm comm into intersection} as $[R_0,\dots, R_n]_S$ and call it a \emph{symmetric commutator subgroup}.

The following surprising result is important for training our models:
\begin{theorem}[\cite{wu2001combinatorial}, Corollary 3.5]\label{thm: partial intersection}
    For any \emph{proper} subset $\{i_0,\dots , i_k\}\subsetneq \{0,\dots , n\}$ the inclusion~\eqref{eq: inclusion of symm comm into intersection} is an equality.
\end{theorem}

This implies that any partial (different from $\cap_{i=0}^n R_i$) intersection of $R_i$ can be easily described and, after defining the suitable distribution, sampled from. In~\cite{wu2001combinatorial} simplicial methods are used to demonstrated that the difference between the symmetric commutator subgroup and the \textbf{full} intersection, as expressed in formula~\eqref{eq: inclusion of symm comm into intersection}, possesses a homotopical nature: 
\begin{theorem}[Wu's formula, \cite{wu2001combinatorial} Theorem 1.7]\label{thm: Wu formula}
    For $n \geq 1$ there is an isomorphism
\end{theorem}
\begin{equation}\label{eq: Wu formula}
    \pi_{n + 1}(S^2) \cong \frac{R_0 \cap R_1 \cap \dots \cap R_n}{[R_0, R_1, \dots, R_n]_S} \cong \frac{R_0 \cap [R_1, \dots, R_n]_S}{[R_0, R_1, \dots, R_n]_S}.
\end{equation}
Although not computational in nature, the formula~\eqref{eq: Wu formula} gives an elegant combinatorial descrtiption of homotopy groups of $S^2$. This description motivates the present work.

\begin{remark}\label{rmk: motivation}
    It is important to note that the problem of classifying whether a given word $w \in \cap_i R_i$ is \emph{trivial} (is an element of the denominator of formula~\eqref{eq: Wu formula}) or \emph{non-trivial} (contributes to a non-zero element of $\pi_* S^2$) is almost as hard as the problem of computing homotopy groups themselves. To our knowledge, there are no known algorithms or approaches to tackle this problem in a general setting, each elements requires an individual method for classification. For instance, for $n = 3, 4$ the group theoretic form of the corresponding generators of the homotopy groups is known~\cite{mikhailov2021homotopy}: 

\begin{center}

    \begin{tabular}{c|c|c}
    $n$&$\pi_{n + 1}(S^2)$ & Non-trivial element \\
    \hline
    2 & $\mathbb{Z} $ & $[x_1, x_2]$ \\
    3 & $\mathbb{Z} / 2$ & $[[x_1, x_2], [x_1, x_2x_3]]$ \\
    4 & $\mathbb{Z} / 2$ & $[[[x_1, x_2], [x_1, x_2x_3]], [[x_1, x_2], [x_1, x_2x_3x_4]]]$ \\
    \end{tabular}
    \captionof{table}{Elements from $\cap_{i=0}^n R_i$, that are not contained within $[R_0,\dots, R_n]_S$}{\label{table: non-trivial elements}}
\end{center}
In our research, we aim to generate elements from the numerator of Wu's formula without providing any information about its denominator to the model, thereby enforcing it not to ``learn'' the pattern of the trivial words in the denominator. Furthermore, during the training phase of our algorithms, we do not intend to utilize the knowledge of non-trivial elements listed in Table~\ref{table: non-trivial elements}, since we want our approaches to be scalable. However, we can consider utilizing the non-trivial elements and the denominator in the evaluation phase to assess the performance of our approaches.
\end{remark}

\section{Generating elements from the intersection of normal closures}
\label{sec: Framework description}

\subsection{Creating datasets from normal closures and their commutator subgroups}\label{subsec: sampling from normal closures}

Generating a consistent synthetic dataset is of paramount importance for data-driven projects. In this regard, one of the biggest challenges is to ensure that the dataset is not only large enough but also exhibits a reliable structure. In our project, we encountered a similar challenge while constructing a training dataset comprising words from $R_i$ and their partial intersections.

In order to effectively train our models, we need to have a certain amount of control on properties of training dataset, for example being able to control to some extent the distribution of lengths of the words and their internal structure.

When it comes to generating words from $R_i$, we discovered that using a \emph{naive} sampling approach, where each $y_k$ in formula~\eqref{eq: definition of R_i} is sampled independently, yielded less desirable properties in the generated words: the adjacent letters of different conjugators $y_k, \ y_{k+1}$ do not interact with each other, resulting in reduced variability within the generated dataset. 

To address this issue, we developed an alternative description of $R_i$, based on balanced bracket sequences and Dyck-$n$ languages. Specifically, to generate a word from $R_i$, we can uniformly sample a balanced bracket sequence~\cite{atkinson1992generating} and then replace every matching pairs of brackets with words $a$ and $b$ such that the word $ab$ is either identity or the rotation of the normal closure's generator. For instance, if $n = 3$, then in case of $R_0$ the opening bracket $a$ could be $x_2$ and the corresponding closing bracket will be $b = x_3x_1$. We refer to this sampling method as \emph{bracket-style sampling}.

We use the notion of Dyck paths to analyse the difference between naive and bracket-based sampling. By \emph{Dyck-$n$ word} we will mean a balanced bracket sequence with brackets of $n$ types. For Dyck-$n$ word $b$ the corresponding \emph{Dyck path} $p(b)$ is a plot that depicts the variation in the number of unbalanced brackets as a function of the position in the word. By \emph{valley} on Dyck path $p(b)$ we will mean a local minimum which is greater than zero.

\begin{figure*}[t]
     \centering
     \begin{subfigure}[b]{0.32\textwidth}
         \centering
         \includegraphics[trim={0.25cm 0.3cm 0.25cm 0.3cm},clip, width = \textwidth]{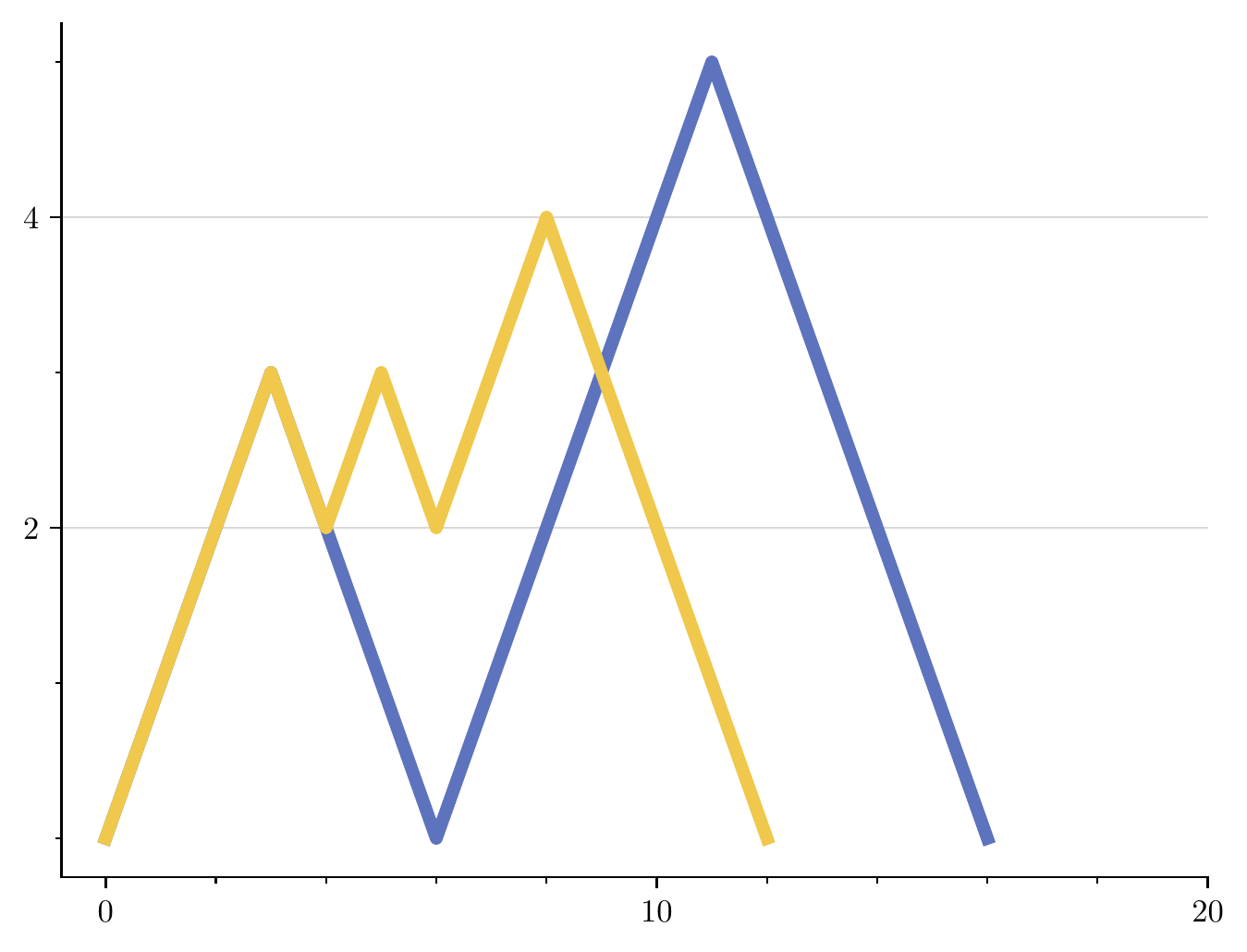}
         \caption{Dyck paths of a typical word\protect\footnotemark}
         \label{fig:paths}
     \end{subfigure}
     \begin{subfigure}[b]{0.32\textwidth}
         \centering
         \includegraphics[trim={0.25cm 0.3cm 0.25cm 0.3cm},clip, width=\textwidth]{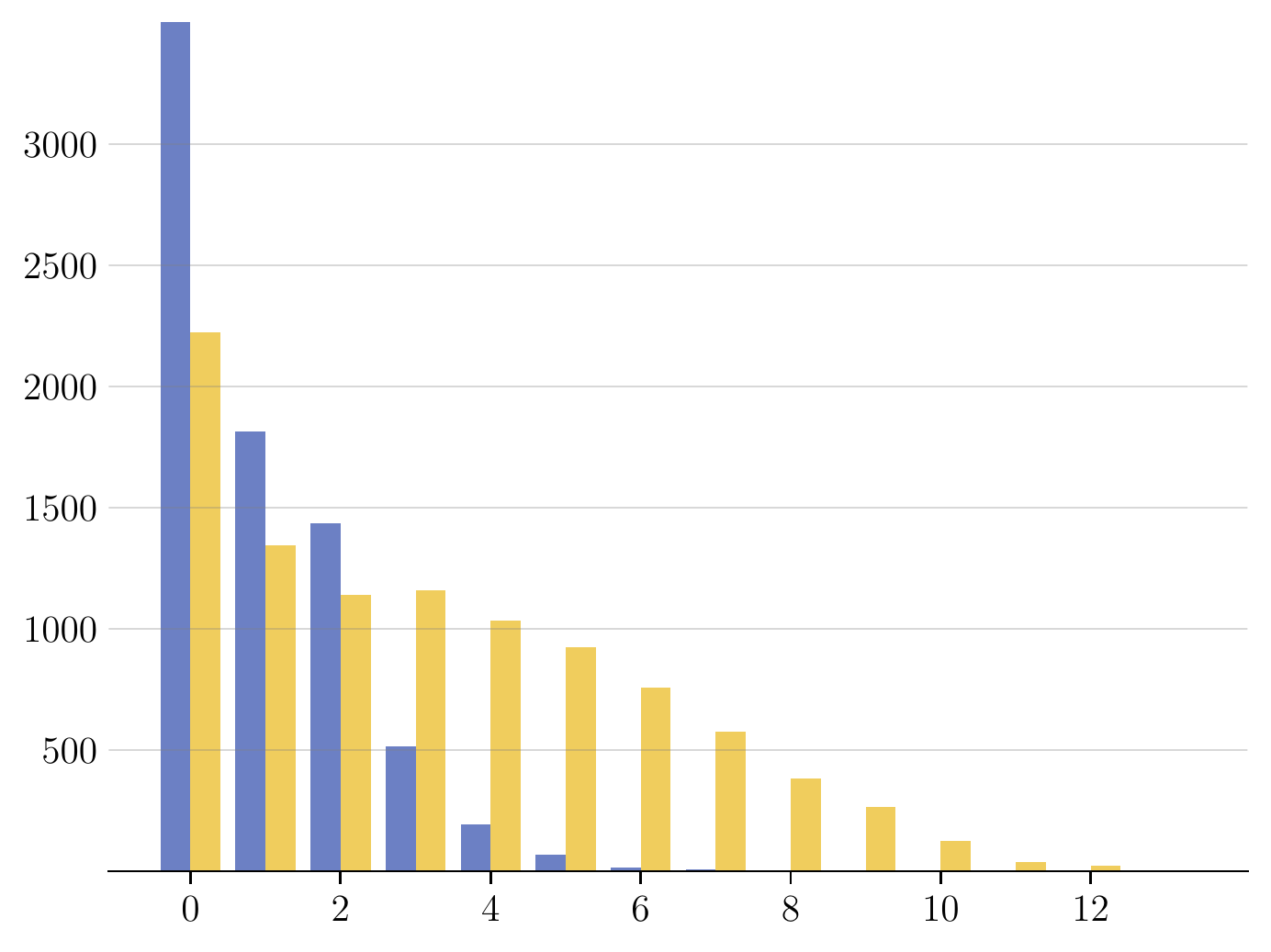}
         \caption{Distribution of valleys}
         \label{fig:valleys}
     \end{subfigure}
     \begin{subfigure}[b]{0.32\textwidth}
         \centering
         \includegraphics[trim={0.25cm 0.3cm 0.25cm 0.3cm},clip, width=\textwidth]{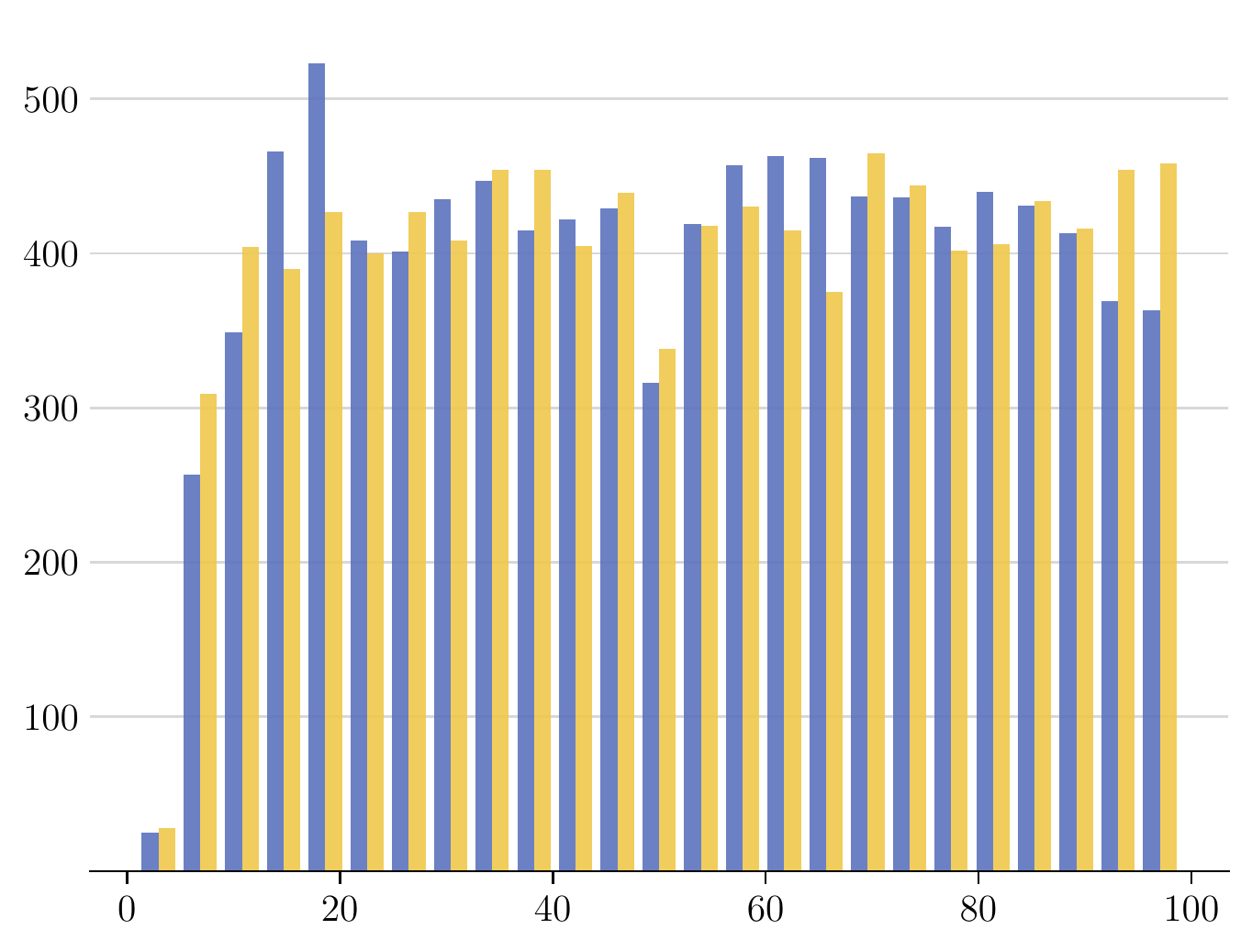}
         \caption{Distribution of lengths}
         \label{fig:lengths}
     \end{subfigure}
        \caption{Comparison of naive (blue) and bracket-style (yellow) sampling. Sample of $10^4$ elements of $R_1$ with $n=3$ with a maximum length of $100$.}
        \label{fig: dyck paths}
\end{figure*}

Figure~\ref{fig: dyck paths} illustrates the difference in distribution of valleys for naive and bracket-style sampling methods. While the distributions of word lengths remain similar, the bracket-style sampling method shows more variety in possible number of valleys in corresponding Dyck paths. Note that each valley corresponds to a mutual cancellation of letters in adjacent conjugators $y_k$ and $y_{k+1}$.

Now, to generate elements from iterative commutator subgroups like $[R_{i_0},\dots, R_{i_k}]$ we use the following procedure. First, we sample a collection of random binary trees that represent products of various commutators. Then for each leaf $j$ of every tree the element of $R_{\sigma(j)}$ is sampled using bracket-style sampling, where $\sigma \in \Sigma_k$ is a random permutation. Note that although commutators of arbitrary complexity can be expressed as product of iterated commutators of the form~\eqref{eq: iterated commutator notation} using Hall-Witt identity (see Appendix~\ref{apdx: background on group theory} for details), we found that sampling not only iterated, but arbitrary commutators, such as $[[x_{i_1},x_{i_2}],[x_{i_3},x_{i_4}]]$ yields better results.

\footnotetext{The words are: $y^{2}[x, z^{-1}]y^{-2}[z^{-1},p][x, z^{-1}][p, z^{-1}]$ with 0 valleys and $p^{2}[z, x][pxp,y]p^{-2}$ with 2 valleys}

Finally, to ensure that the length distribution of words is equal among various partial intersections, we employed the Kolmogorov-Smirnov test. Further details regarding supplementary properties of the sampling methods can be found in Appendix~\ref{apdx: additional experiments and examples}.

\subsection{Baselines}\label{subsec: baselines}
We choose a range of non-neural based algorithms to approach a problem of generating elements from the intersection $\cap_{i = 0}^n R_i$, given algorithms for generating elements from each individual $R_i$.

\paragraph{Random search}
Guided by Theorem~\ref{thm: partial intersection}, the most straightforward and readily implementable method is to randomly sample words belonging to $[R_0,\dots, \hat R_i, \dots, R_n]_S$ and check whenever such these words lie in $R_i$. Notation $\hat R_i$ here means that $R_i$ is excluded from the symmetric commutator subgroup. We are performing this procedure for all $i \neq 0$ simultaneously. Already in random search the difference between naive sampling and bracket-style sampling from $R_i$, discussed in Section~\ref{subsec: sampling from normal closures}, can be seen in terms of number of generated words from intersection per batch (see Appendix~\ref{apdx: additional experiments and examples} for additional experiments).

\paragraph{Evolutionary method}\label{subsubsec: evolutionary method}

As described in Section~\ref{sec: group theory}, the words in $\cap_i R_i$ have are characterized by the following property: replacing any of $x_i$ with $1$ will reduce the word to $1$. This criterion underlies the extremal reformulation of the task of generating words from the intersection of $R_i$ and the following algorithm.

For a normal closure $R = \langle v \rangle^F$ of a word $v$ and a word $w$, let the \emph{algebraic distance} $d(w, R)$ be defined by the length of word after throwing away all occurrences of cyclic permutations of $v$ from $w$. Note that $d(w, R) = 0$ iff $w \in F$. This definition resembles algebraic distance~\cite{taubin1991estimation} between a point $x$ and the kernel of a map $f: \mathbb R^d \to \mathbb R$, defined by $d(x, \text{ker} f) = |f(x)|$.

Next, let's introduce the non-negative objective function
\[
d(w) = \sum_{i = 0}^n d(w, R_i)
\]
and consider the corresponding discrete optimisation problem $\min_{w \in F} d(w)$, solution to which will give an element from $\cap_{i=0}^n R_i$. 

To solve this discrete optimisation problem we apply a modification of the $(1+1)$-evolutionary algorithm~\cite{droste2002analysis} with mutation replacing each of randomly chosen letters in the word with (uniformly) random sampled generator or its inverse. It is worth mentioning that one of our metrics (see Section~\ref{subsec: metrics}) is inspired by the notion of distance between the intersection $\cap_i R_i$ and a given word.

\paragraph{Greedy algorithm}\label{paragraph: greedy}

Experiments with language models (Section~\ref{subsec: model description} below) that generate a full word based on a prefix inspired us to find a non-neural based iterative procedure of appending a next token to a given prefix to generate a full word from the intersection. For each normal closure $R_i$, a stack $S_i$ is created from the given prefix. Specifically, $S_i$ is obtained by removing $x_i$, its inverse, and cyclic permutations (in the case where $i=0$) from the prefix.

To determine the next token in the sequence, the algorithm employs a greedy strategy based on the current tops of the stacks. Each token is assigned points based on two mutually exclusive criteria. The first criterion is that the token candidate should decrease the length of a particular stack. For instance, if the top of a stack $S_i$ is $x_j$, then the next token candidate $x_j^{-1}$ would receive points. The second criterion is that the token candidate should not increase the length of any stack. For example, the token $x_i$ does not increase the length of $S_i$ since it is already excluded from it. Tokens that would reduce the length of the given prefix itself are excluded and assigned zero points. The algorithm then selects the next token with the highest score.

For illustrative purposes we included a toy example of application of this method in Appendix~\ref{apdx: additional experiments and examples}

\subsection{Deep models}\label{subsec: model description}
\paragraph{General architecture details}
We now move to the neural-based algorithms for sampling elements from $\cap_i R_i$. We mostly use decoder-only transformers~\cite{radford2018improving}, \cite{vaswani2017attention} as a main architecture for our models. We will briefly review its features relevant to our problem at hand. 

Decoder-only models are often used in tasks where the input sequence is fixed or pre-processed, such as language modelling~\cite{jurafsky2000speech}, text generation~\cite{li2022pretrained}, which aligns with our task of generating words as sequences of tokens that matches a certain target distribution.

We adopt the classical auto-regressive approach, where each token in the output sequence is generated based on the previous tokens and the decoder's internal state. Specifically, at each time step, the decoder takes in the previous token in the output sequence, along with a context vector that summarizes the encoded input sequence, and generates a new token using a softmax~\cite{christopher2006pattern} function.

We also utilize standard modifications of decoder-only transformer models, such as causal masking, which ensures that the model only attends to previous tokens in the output sequence, and positional embeddings, which provide information about the position of each token in the sequence. These modifications help the model generate high-quality and coherent output sequences~\cite{brown2020language}.

To train such a model we use a cross-entropy loss~\cite{bengio2008neural}:

\begin{equation}\label{eq:cross-entropy}
L(y, \hat{y}) = -\sum_{k}\sum_{i}y_{k,i}\log(\hat{y}_{k,i})
\end{equation}
where $y_{k, i}$, is a true probability distribution (of $k$-th token in the sequence having value $i$), $\hat{y}_{k,i} = \mathbf{Pr}_\theta(s_{k} = i \, | \, s_{1}\dots s_{k-1})$ is the corresponding predicted probability distribution (which depends on model parameters $\theta$) , $s_1,\dots, s_t$ is a sequence of tokens from the training dataset.

We use generators $x_k$ and their inverses $x^{-1}_k$ as tokens, together with the special tokens: $\langle \textbf{bos} \rangle$ to identify the beginning of the sequence, $\langle \textbf{eos} \rangle$ to identify the end of the sequence, and $\langle \textbf{pad} \rangle$ to identify \textit{empty} tokens after the end of the sequence (ignored in formula~\eqref{eq:cross-entropy}). The inclusion of the $\langle \textbf{bos} \rangle$ token has been shown to significantly improve the performance of attention-based models in Dyck-$n$ language recognition~\cite{ebrahimi2020selfattention}, see also~\cite{hahn2020theoretical}, \cite{weiss2021thinking}.

Finally, to sample from models we use two configurations: beam search~\cite{tillmann2003word} with the repetition penalty~\cite{keskar2019ctrl} and nucleus sampling~\cite{holtzman2020curious}.

\subsubsection{Methods}

The key ingredient in our deep learning-based generating approaches is multi-label for each word in the training dataset. Multi-label of the word $w$ provides information about the membership of $w$ in each of the subgroups $R_i$: 
\begin{definition}
  For an element $w$ of a free group on $n$ generators, the \emph{multi-label} is a binary array of length $(n + 1)$ such that its $i$-th position decodes whether $w$ contained in $R_i$ or not.  
\end{definition}

We developed a range of approaches which are different from each other in a way how they treat multi-labels.
\paragraph{Ignore}
This procedure involves the deliberate exclusion of multi-label information during both the training and inference phases. This approach is based on the underlying assumption that the model exhibits sufficient capacity to effectively generalize from a dataset containing words associated with various multi-labels. In essence, this can be described as a zero-shot generation, a problem that has been extensively studied in recent years~\cite{wang2019survey}, \cite{cao2020research}.

\paragraph{Mask}

The second approach involves selectively masking specific components of the model in accordance with the corresponding multi-labels. In this method each attention head is assigned a normal closure, and during the training phase we mask the attention heads that do not correspond to normal closures containing the given word by multiplying their output with zero matrix. In the inference phase we do not mask anything to prompt the model to generate from the complete intersection. Using this masking approach we try to help the model to distribute the knowledge about the normal closures across its parts.

\begin{remark}\label{remark: masking ensemble}
   It is worth noting the \emph{ensemble technique}, which is a particular case of the masking approach. This technique involves training of $n + 1$ models on each normal closure. During the inference phase, the given prefix is fed to each model, and next token is sampled from the joint distribution, which is essentially the sum of the outputs distributions 
   
   \[
   f(w) = \text{Softmax}\left(\sum\limits_{i = 0}^{n} M_{i}(w) \times \chi_{\{ w \in R_i \}}\right),
   \] 
   
   with $\chi$ being a characteristic function of a subset. This approach is functionally equivalent to creating a model that comprises $n + 1$ transformers, with the logits of the last output layers being summed together, as Figure~\ref{fig:masking and ensemble} illustrates. During training, the outputs of each model in ensemble are correspondingly masked according to the multi-label. 

\begin{figure}
    \centering
    \includegraphics[width=0.9\textwidth]{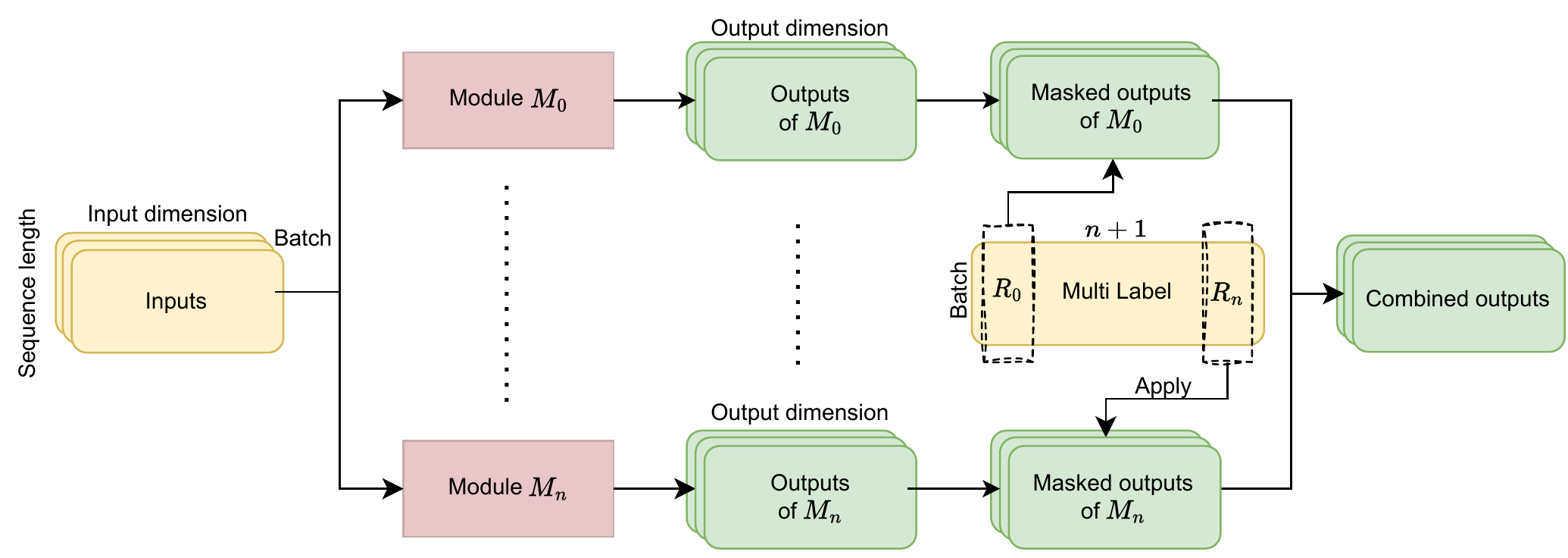}
    \caption{Architecture for a generalized masking model. If each module $M_i$ represents an attention head (or a group of attention heads), then we get a masking model. If the modules $M_i$ represent the complete transformers, we obtain an ensemble model. In both cases the combined outputs is a result of application of the linear transformation (trainable or constant) to the stacked masked outputs.}
    \label{fig:masking and ensemble}
\end{figure}
\end{remark}

\paragraph{Prompt}
Our third approach utilizes multi-label information as a \emph{prompt} for the model. To facilitate this, we expand the vocabulary by introducing tokens for each normal closure $R_0, \dots, R_n$ along with a token $\langle\mathbf{:}\rangle$ that serves as an identifier for the end of the prompt. Prior to inputting the words to the model, we prefix the corresponding prompt to every word. For instance, in the case of a word $x_1x_2\dots x_t$ with multi-label $\texttt{[1, 1, 0, 0]}$, we provide the model $R_0R_1:x_1x_2\dots x_t$ as input. In the evaluation phase we prompt the model with $R_0\dots R_n:$, indicating our intention to obtain a word from the complete intersection. Prompting has become increasingly popular in recent years~\cite{liu2021pretrain} as a means of fine-tuning pre-trained models for specific tasks such as question-answering, summarization, and sentiment analysis. By providing targeted prompts, the model's output can be optimized for specific use cases, making it more accurate and efficient at performing specific tasks. 

\paragraph{Negative Baseline}\label{paragraph: negative baseline}

We can collect the dataset of words from the $[R_0, \dots, R_n]_S$, i.e. a denominator of Wu's formula~\eqref{eq: Wu formula} and train a model to generate elements from it. Although we can not use this method in our comparison, since it was trained only on trivial words (see terminology of Remark~\ref{rmk: motivation}), but we can employ its metrics in order to gain new insights about our primary methods. Furthermore,
negative baseline can be used for so-called \emph{negative knowledge distillation}~\cite{gou2021knowledge}: fine-tuning one of the primary models to predict the outputs that are \textbf{not} predicted by the negative baseline.


\section{Evaluation}\label{sec: evaluation}
Before discussing the results of the evaluation of the proposed models, we explain the relevant metrics and datasets, with other technical implementations details given in Appendix~\ref{apdx: implementation details}.
 
\paragraph{Metrics}\label{subsec: metrics}

There are families of metrics commonly used in natural language processing, which, following~\cite{celikyilmaz2020evaluation}, may be divided in three groups: human-centric evaluation metric, automatic metrics and machine-learned metrics. We decide to concentrate our attention on automatic metrics, since our primary objective is to understand the quantitative performance of neural-based algorithms in a task of generating the words from $\cap_i R_i$.

Our main metric, the \emph{completion ratio}, is an average number of the successfully generated words per batch. Specifically, for a batch of randomly generated words prefixes, the model generates a set of possible outputs based on each prefix. The number of generated outputs that simultaneously belong to all normal closures $R_i$ is then computed, and this number is divided by the size of the batch. A higher completion ratio signifies a greater proficiency of the model in generating desired outputs.

In addition to CE loss~\eqref{eq:cross-entropy}, we also utilize a \emph{reduction ratio} for guiding the training process. This metric measures the distance (see Section~\ref{subsubsec: evolutionary method}) from the generated output to the given $R_i$ divided by the length of the input, batch averaged and averaged over all $i$.

Reduction ratio sits ``in between'' the loss and the final completion ratio metric: once the loss on the validation dataset has stabilized, and the completion ratio remains relatively low, the changing rate of the reduction ratio becomes informative. This changing rate indicates how the model's output gradually approaches the numerator of Wu's formula (see Figure~\ref{fig: Metrics}).

\begin{remark}
We also notice that although for future projects of generation \emph{non-trivial} words human evaluation is still out of scope due to our constrained knowledge of the intricate structure of the homotopy groups. While there is a relative limited number of mathematicians deeply immersed in this topic, the machine-learned metrics (for instance, probability of the word being contained in the symmetric commutator subgroup) may be beneficial and play significant role. 
\end{remark}

\paragraph{Datasets}\label{paragraph: datasets}
In order to maximize the performance of our models, specifically in terms of the number of generated words from random prefixes, we adopt an approach where overfitting is not a primary concern. To achieve this, we employ an infinite, online generated training dataset that encompasses words from $[R_0,\dots,\hat R_i,\dots, R_n]_S$ for all $i$, i.e. the multi-labels associated with these words contain only a single zero. The validation dataset is also generated in an online fashion, drawing samples from the same distribution. For a negative baseline model, described in Section~\ref{paragraph: negative baseline}, both training and validation dataset consist of words from the full symmetric commutator subgroup $[R_0,\dots, R_n]_S$.

For an auxilary tests, we have assembled small datasets consisting of cyclic permutations of generators of $\pi_n S^2$ for $n=3,4$ from Table~\ref{table: non-trivial elements}. By evaluating our metrics on these cyclic permutations, we can assess the models proficiency in generating \emph{non-trivial} words.

\begin{table}[t] 
\centering
\begin{tabular}{c| c | c | c | c | c | c }
    & Random & Evo & Greedy & Ignore & Mask & Prompt \\ 
    \midrule
    $n=3$ & $0.058$ & $0.035$  & 0.479 & $\textbf{0.36}$  & $0.16$  & $0.30$  \\
    $n=4$ & $0.032$ & $0$ & 0.016 & $0.356$  & $0.584$  & $\textbf{0.684}$ \\
    $n=5$ & $0.025$ & $0$ & 0.004 & $0.263$  & $0.487$  & $\textbf{0.550}$ \\
     
\end{tabular}

\caption{Completion ratio of developed algorithms for various number of generators $n$.}
\label{tab:completion-ratio}
\end{table}

\subsection{Results}
\begin{figure*}[t]
     \centering
     \begin{subfigure}[b]{0.32\textwidth}
         \centering
         \includegraphics[trim={0.25cm 0.3cm 0.25cm 0.3cm},clip,width = \textwidth]{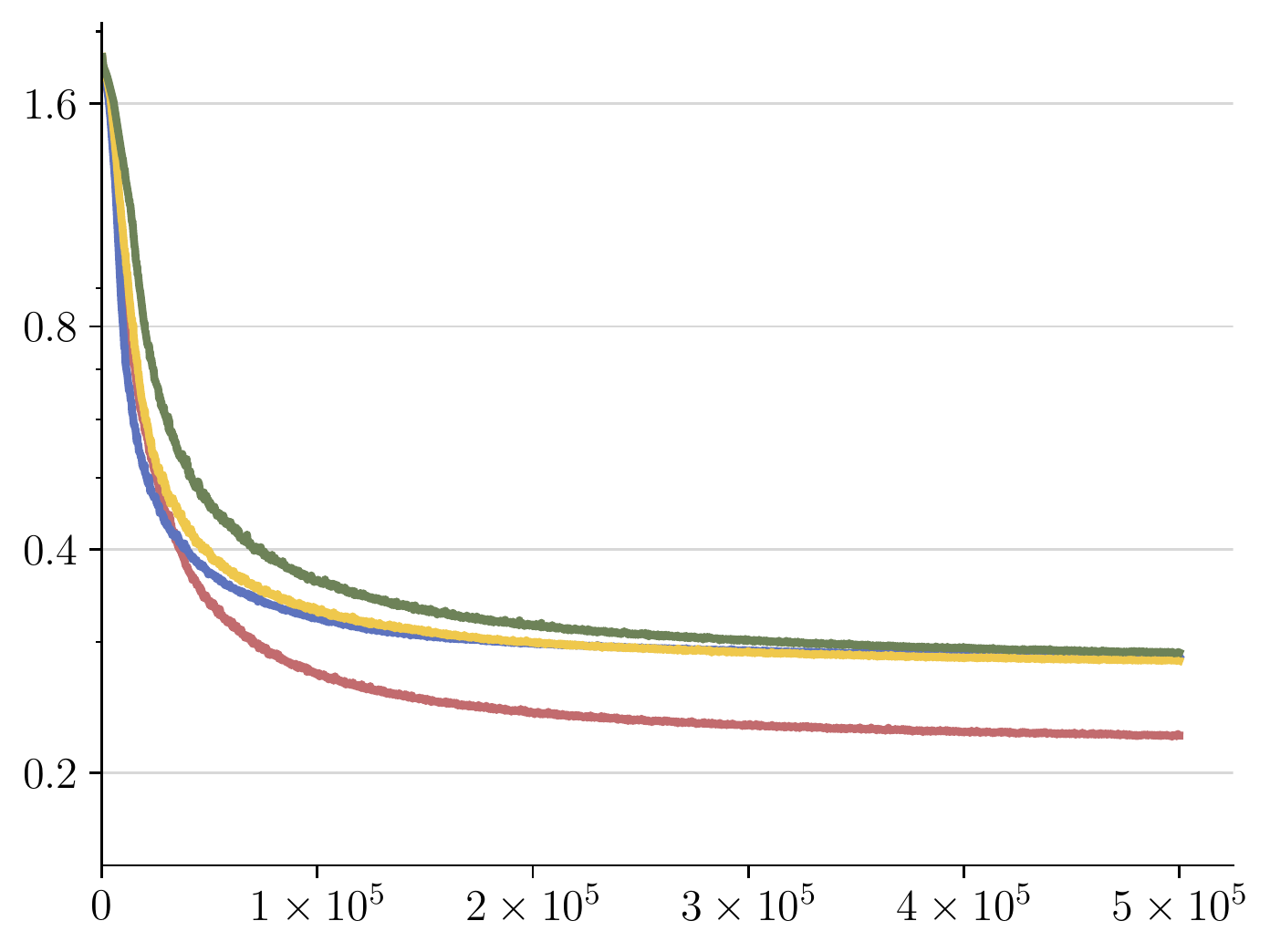}
         \caption{CE loss}
         \label{fig:completion}
     \end{subfigure}
     \hfill
     \begin{subfigure}[b]{0.32\textwidth}
         \centering
         \includegraphics[trim={0.25cm 0.3cm 0.25cm 0.3cm},clip,width=\textwidth]{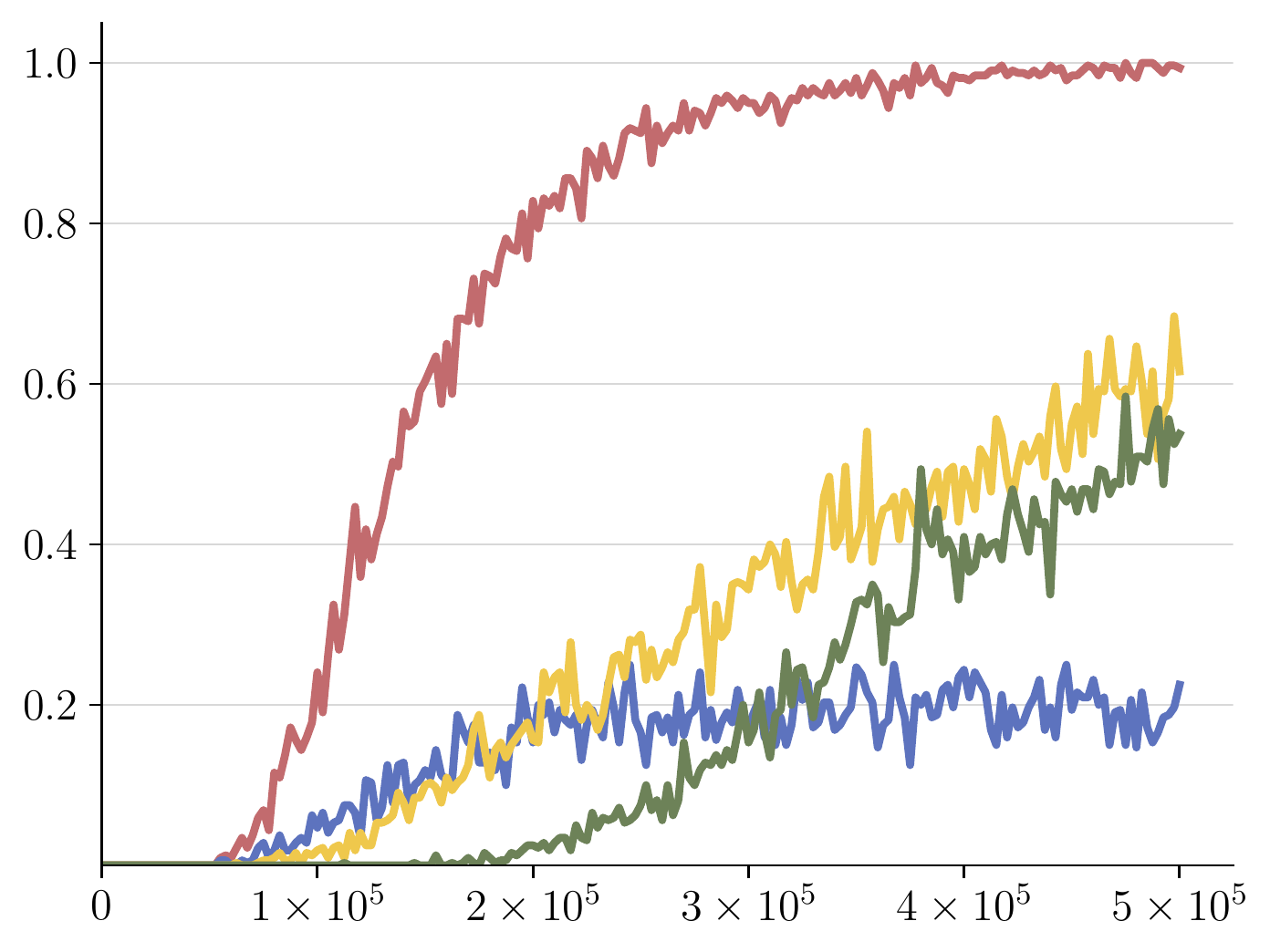}
         \caption{Completion ratio}
         \label{fig:reduction}
     \end{subfigure}
     \hfill
     \begin{subfigure}[b]{0.32\textwidth}
         \centering
         \includegraphics[trim={0.25cm 0.3cm 0.25cm 0.3cm},clip,width=\textwidth]{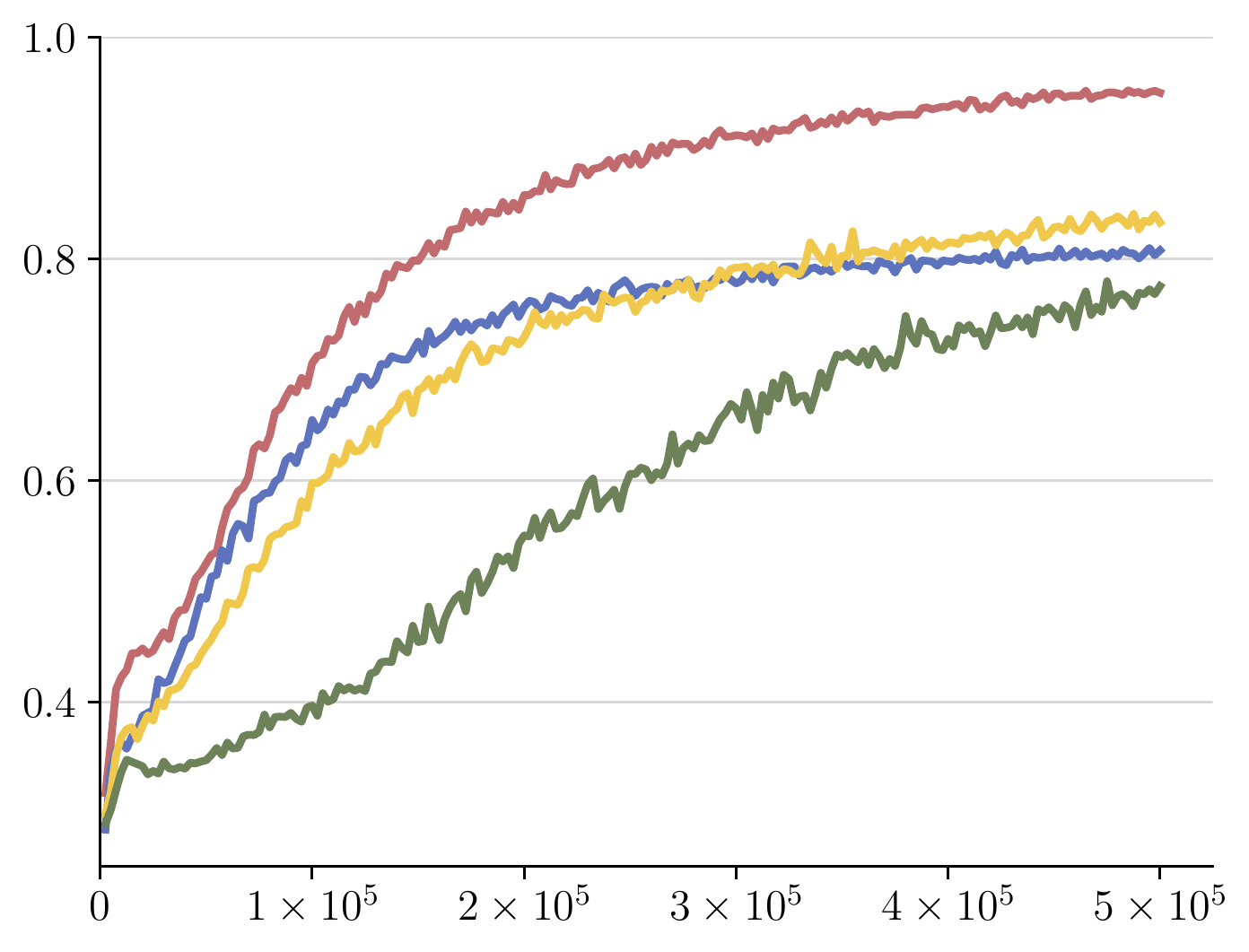}
         \caption{$1-$Reduction ratio}
         \label{fig:non-trivial}
     \end{subfigure}
        \caption{Illustration of the learning process for \emph{negative baseline} (red), \emph{prompt} (yellow), \emph{masking} (green) and \emph{ignore} (blue) methods, $n = 4$.}
        \label{fig: Metrics}

    \begin{subfigure}[t]{0.32\textwidth}
         \centering
         \includegraphics[trim={0.25cm 0.3cm 0.25cm 0.3cm},clip,width = \textwidth]{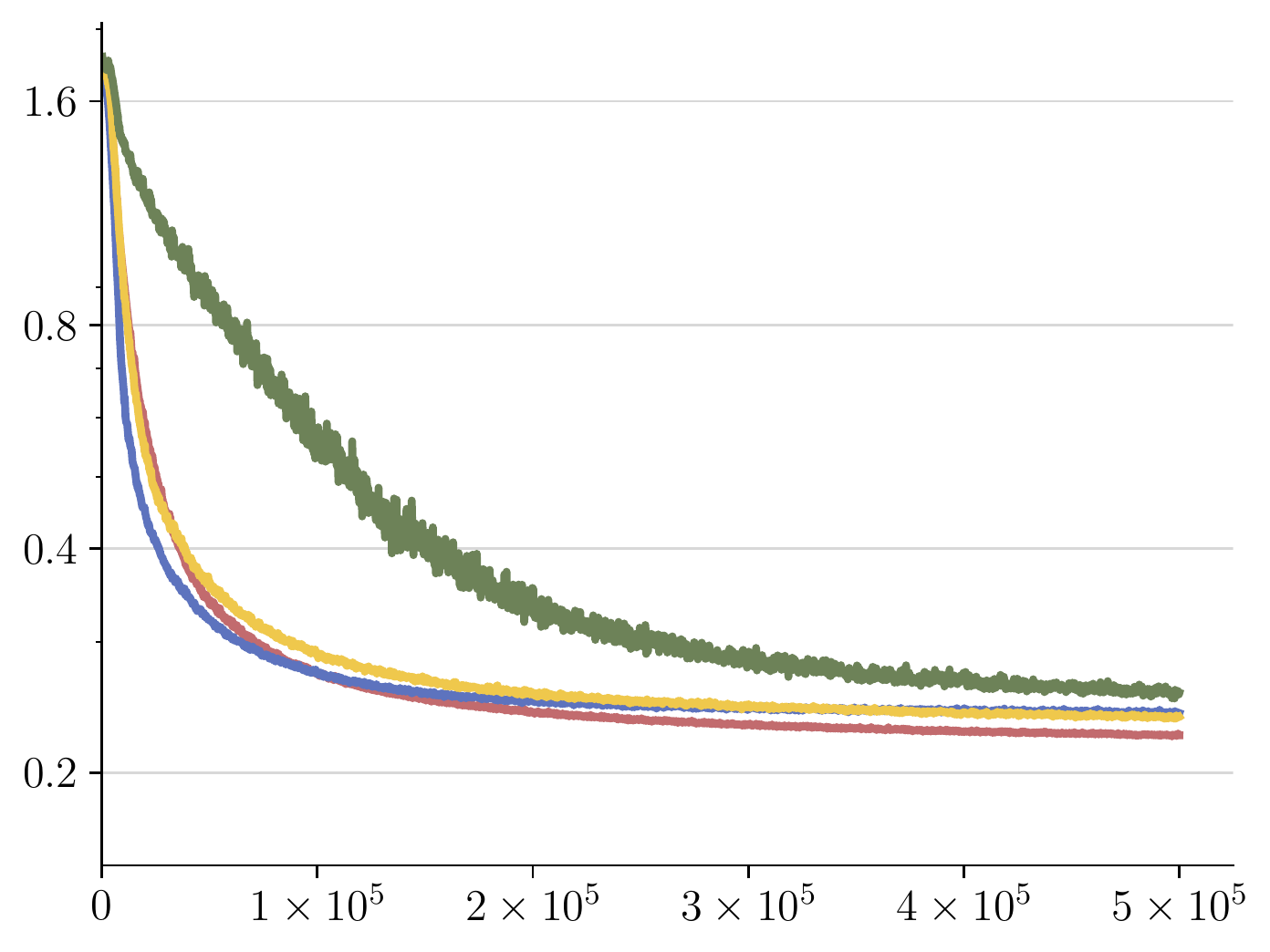}
         \caption{$[R_0,\dots, R_n]_S$}
         \label{fig:triv loss}
     \end{subfigure}
     \hspace{50pt}
     \begin{subfigure}[t]{0.32\textwidth}
         \centering
         \includegraphics[trim={0.25cm 0.3cm 0.25cm 0.3cm},clip,width=\textwidth]{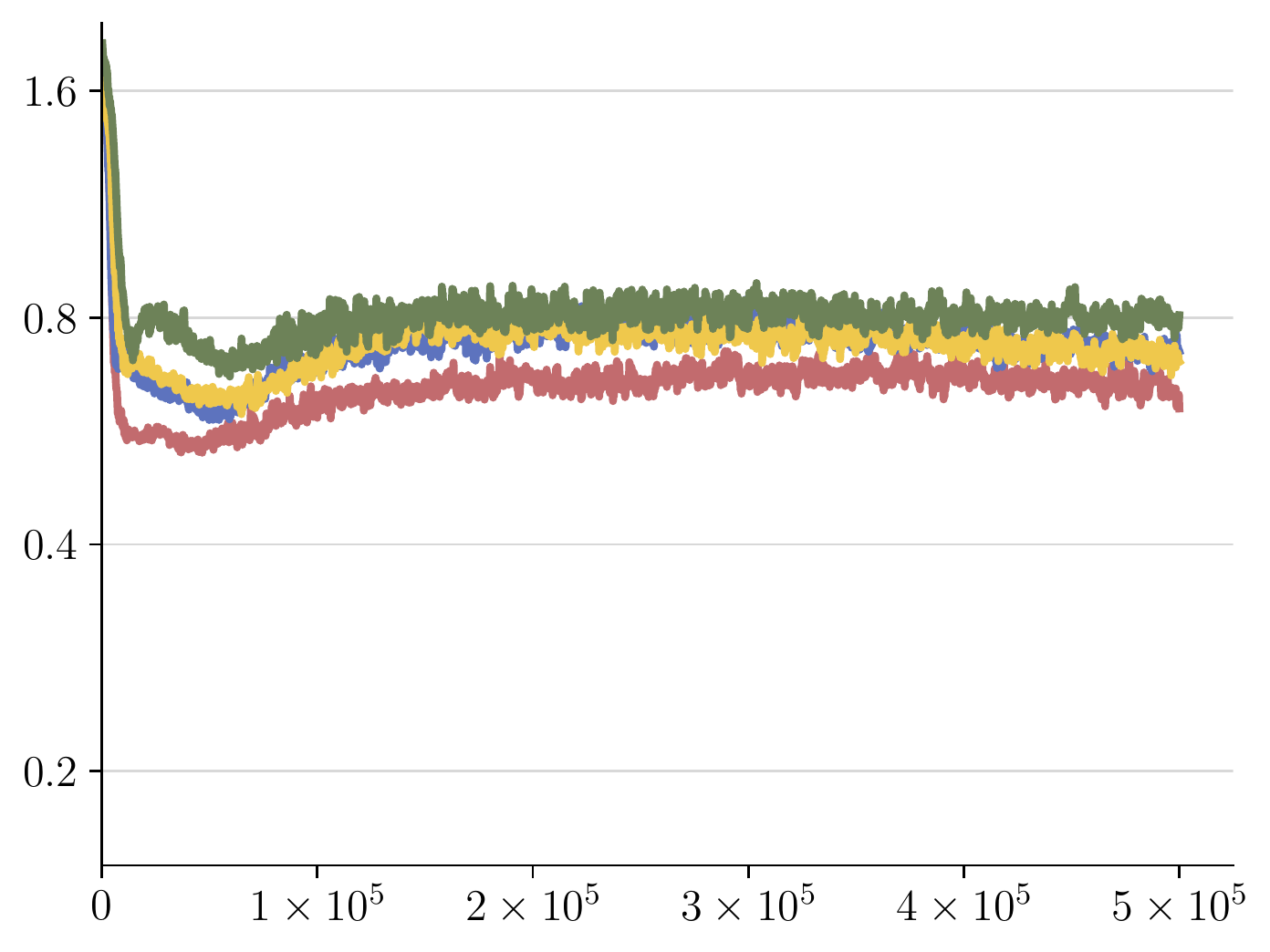}
         \caption{$\bigcap_{i=0}^n R_i \backslash[R_0,\dots, R_n]_S$}
         \label{fig:non-triv loss}
     \end{subfigure}
        \caption{Comparison of CE loss per epoch on validation for the trivial and non-trivial words, $n = 4$.}
        \label{fig: CE loss on triv and non-triv}
\end{figure*}

We compare our deep learning approaches with non-neural-based baselines using completion ratio, see Table~\ref{tab:completion-ratio}. We also report the metric values on validation during training progress (Figure~\ref{fig: Metrics}), to illustrate the difference between various multi-label approaches. We report our results for $n \leq 5$ generators, due to the available computation resources. We also notice (see discussion below), that for a higher number of generators the high length of words from the full intersection $\cap_i R_i$ suggests to use another, more sophisticated representation of words.

\vspace{-5pt}

Table~\ref{tab:completion-ratio} provides compelling evidence of the scalability and superior performance achieved by language models when compared to the baselines. These results signify a promising avenue for future applications of such models in understanding simplicial homotopy groups.  In terms of comparing the performance of the methods between each other, the prompt-based approach outperformed the other methods on the validation dataset. Finally, Figure~\ref{fig: CE loss on triv and non-triv} illustrates that the non-trivial words are significantly more challenging (as expected) for the models to generate.

\section{Conclusion}
We evaluate perspectives of using machine lerning in generating cycles in simplicial groups. From a machine learning perspective, the problem can be reformulated as the task of generating elements from the intersection of subsets of (a form of) Dyck-$n$ words. Some classical algorithms are implemented, some are developed from scrath and used as baselines for deep learning algorithms. Built of the idea of ensemble of generators, the additional multi-label information is added to the training dataset, which allows the single model to work as a generalization of an ensemble. Resulting models, in contrast to baselines, are scalable, and will serve as a building blocks for future algorithms specialized in sampling from homotopy groups of spaces.
\paragraph{Further directions}\label{par: further directions}
Employing a more sophisticated tokenization scheme could lead to a more compact representations of elements in higher simplicial homotopy groups. One of the choices would be tokenizing commutator brackets and delimiters to encode commutator presentation of the words from $\cap_i R_i$. However, in the present work we choose the generators of free group as tokens for its simplicity. Even more compact, but at the same time more elaborate and theoretically involved representation is given by \emph{$\Lambda$-algebra}~\cite{bousfield1966mod}. Our ongoing research is devoted to employ machine learning in building connections between cycles in simplicial groups and elements in $\Lambda$-algebra.

The ultimate objective of our research is to develop a method capable of sampling elements from simplicial homotopy groups. However, it is crucial for such method to distinguish between the intersection of $R_i$ and the symmetric commutator subgroup. We believe that further enhancement of our deep learning techniques will bring this goal closer to fruition.

\medskip

\clearpage
\newpage
\small
\printbibliography

\appendix
\clearpage
\newpage
\section{Background on group theory}\label{apdx: background on group theory}
As a reminder, we will present a series of basic definitions from group theory. These definitions can serve as points of reference from the main body of the paper.

\begin{definition}
A \emph{group} $G$ is a set equipped with a binary operation $\mu: G \times G \to G$, usually denoted just by $(g_1, g_2) \mapsto g_1g_2$, satisfying the following axioms:
\begin{itemize}

\item \textbf{Associativity}: For all $g_1, g_2, g_3 \in G$, $(g_1g_2)g_3 = g_1(g_2g_3)$.
\item \textbf{Identity Element}: There exists an element $1 \in G$, called the \emph{identity element}, such that for all $g \in G$, $g1 = 1g = g$.
\item \textbf{Inverse Element}: For every element $g \in G$, there exists an element $g^{-1} \in G$, called the \emph{inverse element} of $g$, such that $gg^{-1} = g^{-1}g = 1$.
\end{itemize}
\end{definition}

An arbitrary group can be expressed as a quotient (Definition~\ref{def: quotient}) of a free group. 
\begin{definition}
Let $X$ be a set. The \emph{free group} $F$ on $X$ is a group defined as follows: $F$ consists of all reduced words (finite sequences) on elements of $S$ and their formal inverses, equipped with the operation of concatenation followed by reduction. The reduction operation eliminates pairs of the form $xx^{-1}$ and $x^{-1}x$. The identity element of $F$ is the empty word.
\end{definition}
Free groups, in comparison with arbitrary groups, have a simple algebraic structure, so as their \emph{subgroups}, i.e. subsets closed under group operations.
\begin{theorem}[Nielsen–Schreier theorem, \protect{\cite[Section~14.3]{kargapolov1979fundamentals}}]
    Every subgroup of a free group is free.
\end{theorem}
Subgroups $R_i$ of Definition~\ref{eq: definition of R_i}, that we are dealing with, are defined in a way that they are closed under conjugation operation.
\begin{definition}
Let $G$ be a group. A subgroup $N$ of $G$ is called a \emph{normal subgroup} if and only if for every $g \in G$ and $n \in N$, the element $n^g = g^{-1}ng$ is also in $N$.
\end{definition}
\begin{definition}
Let $G$ be a group and $S$ be a subset of $G$. The \emph{normal closure} of $S$ in $G$, denoted by $\langle S \rangle^G$, is the smallest normal subgroup of $G$ that contains $S$. It is the intersection of all normal subgroups of $G$ containing $S$.
\end{definition}
Wu's formula~\eqref{eq: Wu formula} involves a quotient of an intersection of $R_i$'s by their symmetric commutator subgroup.
\begin{definition}\label{def: quotient}
Let $G$ be a group and $N$ be a normal subgroup of $G$. The \emph{quotient group}  of $G$ by $N$, denoted by $G/N$, is the set of \emph{cosets} of $N$ in $G$, i.e. subsets of $G$ of the form $gN = \{gn|n\in N\}$, with group operations induced from $G$.
\end{definition}
Finally, some of our methods rely on the presentations of $R_i$ as kernels of \emph{homomorphisms} $d_i$, i.e. maps between groups that preserve group operations, see also Formula~\eqref{eq: ker d_i = R_i}.
\begin{definition}
    A \emph{kernel} of a group homomorphism $f: G \to H$ is the following subset:
    \[
    \ker f = \{g \in G | f(g) = 1\}.
    \]
\end{definition}
Note that the kernel of a group homomorphism $G\to H$ is always a normal subgroup of $G$.

Throughout the paper, we utilize the notion of commutators and conjugations, which were introduced in Section~\ref{sec: group theory}, to provide a shorter representation of words in a free group. There is a number of \emph{identities}, involving commutatrors, that hold for all elements of any group, for instance:
\begin{gather*}
    [x,zy] = [x,y][x,z][x,z,y] \\
    [xz,y] = [x,y][x,y,z][z,y] \\
    [x,y^{-1},z]^y[y,z^{-1},x]^z[z,x^{-1},y]^x = 1.  
\end{gather*}

The last identity, called the \emph{Hall-Witt identity} is a ``non-abelian'' version of the Jacobi identity in Lie algebras. As mentioned at the end of Section~\ref{subsec: sampling from normal closures}, these identities enable us to express every commutator as a product of iterated commutators \textbf{of the same length}. For example, for commutator  $[[a,b],[c,d]]$, the Hall-Witt identity reads:
\[
[a,b,[c,d]]^{b^{-1}}[b^{-1},[d,c],a]^{[c,d]}[c,d,a^{-1},b^{-1}]^a = 1
\]
and since
\[
[b^{-1},[d,c],a] = [d,c,b^{-1},a^{[b^{-1},[d,c]]}]^{-1},
\]
we get
\[
[[a,b],[c,d]] = [c,d,a^{-1},b^{-1}]^{-ab}[d,c,b^{-1},a^{[b^{-1},[d,c]]}]^{[c,d]b}.
\]

\section{Background on simplicial homotopy theory}\label{apdx: background on simplicial homotopy theory}
In this section we will provide a brief overview of the historical development of the concept of homotopy groups (see also~\cite{hilton1988brief}). Subsequently, we will delve into the basics of simplicial homotopy theory and the description of Wu's formula within the context of simplicial groups. 
\subsection{History overview}
As a distinguished area of mathematics, algebraic topology was surfaced since the end of XIX century. During this time, mathematicians realized that certain important properties of geometric objects, 
like properties of vector fields on a manifolds or solutions of differential equations on them are controlled by \emph{topological} properties of the underlying objects of study, rather than geometric or analytic characteristics. This realization prompted the need for algebraic invariants that could discern between classes of spaces with similar topological properties and aid in their classification. Historically, one of the first invariant was \emph{homology} $H_*$ of space, shortly followed by a fundamental group $\pi_1$.

Such invariants were able to distinguish spaces \emph{up to homotopy equivalence}. Specifically, two maps $f,g: X\to Y $ are called homotopic, if there is a family of maps $h_t: X \to Y, \ t\in [0,1]$, continuous in $t$, which ``connects'' $f$ to $g$ in a sense that $h_0 = f, \ h_1 = g$. Similarly, two spaces $X,Y$ are \emph{homotopy equivalent} if there exist maps $X\leftrightarrows Y$, which compositions with each other are homotopic to identities.

Generalizing $\pi_1$, the series of higher invariants $\pi_*$, called the \emph{homotopy groups}, were defined as homotopy classes of maps from higher dimensional spheres $S^n$ to a space of interest. Stronger than homology, homotopy groups capture the substantial part of a homotopy type of space, which we nowadays call a \emph{weak homotopy type}. Extremely difficult to compute and comprehend, till today, homotopy groups remain the main object of interest in algebraic topology. Various flavors of homotopy theory have been developed, including stable homotopy, $v_n$-periodic homotopy and others, each with corresponding homotopy groups that are more tractable and accessible than the classical ones.

In parallel with the development of various algebraic invariants for classifying spaces, the concept of space itself has undergone several formalizations within homotopy theory. Initially, the notion of a space referred to a \emph{topological space}, which consisted of a set equipped with a chosen collection of subsets called a \emph{topology}. Later the concept of an (abstract) simplicial complex emerged as a combinatorial notion of space, suitable for defining homology groups. This notion was subsequently refined to that of a simplicial set.

Homotopy theory of simplicial sets, which is both combinatorial and algebraic in nature, now serves as the foundation for more advanced concepts in abstract homotopy theory, such as simplicial model categories, quasicategories, simplicial localizations, and others. Within the scope of the current paper, a specific type of simplicial sets called \emph{simplicial groups} holds particular significance. Simplicial groups serve as a crucial link between simplicial homotopy theory and group theory, offering a concrete, algebraic, and computational-friendly framework.

\subsection{Simplicial groups and Wu's formula}
We will now describe the relevant part of simplicial homotopy theory in greater detail. Simplicial set (formally described as a functor $\Delta^{op} \to \Set$) can be visualized~\cite{mac2013categories}[Ch. 7.5] as a diagram of sets
\[
X : \hspace{1cm}
\begin{tikzcd} {\pi_0 X} & \arrow{l}[swap]{d_0} 
X_0  \arrow{r}{s_0} 
& X_1  \arrow[shift left=4mm]{l}[swap]{d_1} \arrow[shift left=-4mm]{l}[swap]{d_0} 
\arrow[shift left=4mm]{r}{s_0} \arrow[shift left=-4mm]{r}{s_1}
& X_2   \arrow[shift left=8mm]{l}[swap]{d_2} \arrow{l}[swap]{d_1} \arrow[shift left=-8mm]{l}[swap]{d_0} 
& \dots,
\end{tikzcd}
\]
with maps $d_i$, called \emph{faces} and $s_j$, called \emph{degeneracies} which satisfy the \emph{simplicial identities}:
\begin{align}\label{eq: simplicial relations}
\begin{split}
    d_id_j = d_{j-1}d_i , \ i < j \\
    s_is_j = s_js_{i-1} , \ i > j \\
    d_is_j = 
    \begin{cases}
    s_{j-1}d_i, \ i<j \\
    \mathrm{id}, \ i\in \{j,j+1\} \\
    s_jd_{i-1}, \ i > j + 1
    \end{cases}.
\end{split}
\end{align}
Drawing intuition from simplicial complexes, each $X_n$ can be thought of as a set of $n$-simplices, with maps $d_i: X_n \to X_{n-1}$ indicate an $i$-th face of an $n$-simplex and maps $s_j: X_{n-1} \to X_n$ indicate ways how to consider $(n-1)$-simplex as a (degenerate) $n$-simplex.

Simplicial sets were introduced by Eilenberg and Zilber~\cite{eilenberg1950semi} (originally under the name of semi-simplicial complexes) for the needs of formalization of homology theory. When the abstract homotopy theory became more articulated in the fundamental works of Gabriel and Zisman~\cite{gabriel2012calculus} and Quillen~\cite{quillen2006homotopical}, it was shown that the corresponding homotopy theory is equivalent to the classical homotopy theory of topological spaces. Although for some concrete work in homotopy theory simplicial sets may better suited than topological spaces/cell complexes, regarding computations of homotopy groups, bare simplicial sets have some difficulties when it comes to computation of homotopy groups. For instance, the naive notion of homotopy is not an equivalence relation for all simplicial sets.

One of the possible solutions is to consider \emph{simplicial groups}, i.e. simplicial sets which have a group structure on the sets of their simplices, such that the maps $d_i, s_j$ preserve this structure. A notable advantage of simplicial groups is the very explicit description of their homotopy groups: if $G$ is a simplicial group, then
\begin{equation}\label{eq: simplicial group homotopy groups}
  \pi_n G = \frac{\cap_{i = 0}^n \ker d_i}{\mathrm{im}\, \partial_{n+1}},  
\end{equation}

where $\partial_{n+1}$ is a restriction of $d_{n+1}$ to the intersection of kernels of the next dimension. So the homotopy groups of $G$ are identified with the quotients of intersections of kernels of face maps $d_i$ by the images of the last face maps (restricted to a suitable intersection).

To get the Wu's formula for $\pi_n S^2$ from~\eqref{eq: simplicial group homotopy groups}, we will consider a particular simlicial group (more precisely, the functor $s\Set \to s\Gr$) called the \emph{Milnor construction}~\cite{curtis1971simplicial}[Ch. 4]. Let $X$ be a connected simplicial set with a basepoint $*$. The Milnor construction for $X$, denoted by $F[X]$ is a simplicial group obtained by level-by-level-wise application of the free group functor to $X$, quotient by the relation $* = 1$. This means that $n$-simplices in Milnor's construction $F[X]$ form free group on $n$-simplices $X_n$ with the identity of the group identified with a basepoint $*$. The reason for introducing this construction is that for a connected simplicial set $X$ the Milnor's construction $F[X]$ has a homotopy type of loops over the suspension of $X$. For $X = S^1$ on the level of homotopy groups it means that
\[
\pi_n F[S^1] = \pi_{n+1} S^2.
\]

The simplicial circle $S^1$ consists of one non-degenerate $0$-simplex, which is a basepoint $*$ and one non-degenerate $1$-simplex $\sigma$ (plus the degeneracy $s_0(*)$), and there are no non-degenerate simplices in all higher dimensions. 

\[
\begin{tikzcd}[column sep=large]
    * 
    \arrow["{\sigma}","\mathstrut"{name=0, anchor=center, inner sep=0}, from=1-1, to=1-1,
        to path={let \p1=($(\tikztostart.north)-(\tikztostart.south)$),\n1={scalar(asin(0.5*\y1/1.5em))} in
            (\tikztostart.south) arc[start angle=-\n1,end angle=-360+\n1,radius=1.5em]\tikztonodes (\tikztotarget)}]
\end{tikzcd}
\]

Now simplicial identities~\eqref{eq: simplicial relations} show that 
\[
     F[S^1]_n = F[\{x_0, \dots, x_{n-1}\}]
\]
and faces are given by
\begin{equation*}
    d_0(x_j)=
    \begin{cases}
    1, \ j = 0 \\
    x_j, \ j > 0
    \end{cases}\ 
    d_i(x_j) = 
    \begin{cases}
    x_j, \ j < i \\
    x_{i - 1}x_i, \ j \in \{i - 1, i\} \\
    x_{j - 1}, \ j > i
    \end{cases}\ 
    d_n(x_j)=
    \begin{cases}
    x_j, \ j < n - 1 \\
    1, \ j = n - 1
    \end{cases}.
\end{equation*}

It is shown in~\cite{wu2001combinatorial} that the kernels $\ker d_i$, after an appropriate change of basis and relabeling, can be identified with normal subgroups
\begin{align}\label{eq: ker d_i = R_i}
     \ker d_0 = R_0 = \langle x_1\dots x_n \rangle^F
     \ker d_i = R_i = \langle x_i \rangle^F, i > 0 ,
\end{align}

and the image in the denominator of Formula~\eqref{eq: simplicial group homotopy groups} with their \emph{symmetric commutator subgroup}

\[
     \mathrm{im}\, \partial_{n + 1} = [R_0,\dots, R_n]_S = \prod_{\sigma \in \Sigma_{n+1}}[R_{\sigma(0)}, \dots, R_{\sigma(n)}].
\]

This is the idea behind the Wu's formula.

\section{Implementation details}\label{apdx: implementation details}
To operate with words in free groups we implemented a small \texttt{Python} module \texttt{freegroup}\footnote{Code available at \url{https://github.com/ml-in-algebraic-topology/freegroup}}. In this module we use lists of integers as the representation of words with the following mapping: $x_k \mapsto k, x^{-1}_k \mapsto -k$.

As a base for our models we used GPT2~\cite{radford2019language} implemented by \texttt{huggingface}~\cite{wolf2020transformers}. For $n = 3$ we used the model with $\sim 2\cdot 10^6$ parameters and for $n = 4, 5$ we used $\sim 50 \cdot 10^6$ parameters. We employed standard training framework with iterating through available data and gathering examples to batches. We trained models for $2 \cdot 10^5, 5 \cdot 10^5$ iterations with batches of sizes $64, 16$ for $n = 3$ and $n = 4, 5$ cases, respectively. For the optimization we used AdamW~\cite{loshchilov2019decoupled} implemented by \texttt{PyTorch}~\cite{paszke2017automatic} with learning rates: $10^{-4}$ for $n 
 = 3$ and $10^{-6}$ for $n = 4, 5$. For the inference we used generation algorithms implemented by \texttt{huggingface}. During sampling $\texttt{top\_p} = 0.9$; during beam search $\texttt{num\_beams} = 5$, $\texttt{repetition\_penalty} = 1.2$.

Dataset parameters are the following: 
\begin{itemize}
    \item maximal length of a word = $200, 400, 600$,
    \item maximal lengths of a word from $R_0$ are $10, 9, 8$ for $n = 3, 4, 5$ respectively,
    \item  maximal lengths of a word from $R_i, \ i > 0$ are $30, 30, 30$ for $n = 3, 4, 5$ respectively.
\end{itemize}

\section{Additional experiments, models and examples}\label{apdx: additional experiments and examples}
\subsection{Other models}
In addition to the deep learning methods discussed in the main body of the paper, we explored several other trainable and non-trainable machine learning approaches for word generation.

\paragraph{Fixed embeddings} 

Alongside the development of deep learning models, which are in fact suitable deep feature extractors coupled with a classifiers, we also investigated the usefulness of non-trainable embeddings and experimented with different optimizers and machine learning classifiers applied on top of them. The underlying motivation for this investigation was to transform the problem setting from a discrete one to a continuous one, as most optimization and machine learning methods operate more effectively in continuous domains.

The variant of the non-trainable embeddings that we investigated is based on the matrix representation of free group. Using the group homomorphism $\rho: F_2\to SL_2(\mathbb R)$, called \emph{Sanov representation}~\cite{zubkov1998matrix}
\begin{equation}\label{eq:Sanov_repr}
    x\mapsto \begin{pmatrix}1 & 2 \\ 0 & 1\end{pmatrix}, \ y\mapsto \begin{pmatrix}1 & 0 \\ 2 & 1\end{pmatrix},
\end{equation}
the free group $F_2$ on two generators $x$ and $y$ can be embedded in a space of two-by-two matrices with real coefficients (in fact, into a subgroup of matrices $M$ with $\mathrm{det} M = 1$) 

The embedding of a free group $F$ on an arbitary number of generators $n$ is obtained by inclusion
\[
F \to F_2, \ x_i \mapsto [x,y^i].
\]
\begin{remark}
Note that for a word $w \in F$ the magnitudes of coefficients in the matrix $\rho(w)$ depends exponentially on the length of $w$, making this embedding practically unfeasible. Nevertheless we tested an optimization on this embedding space to evaluate a potential of non-trainable embeddings on a small scale.
\end{remark}

One of the advantages of matrix embeddings is that the homomorphisms $d_i: F \to F$, such that $\ker d_i = R_i$, as in Formula~\eqref{eq: ker d_i = R_i} can be reinterpreted as a matrix operation of nullifying $i$-th column denoted $f_i$, $i>0$. At the same time, since $\rho$ preserves the group operation, verifying whether $w \in R_i$ is equivalent to checking if $f_i \rho(w)$ is equal to the identity matrix $I$.

We investigated various optimization procedures on this feature space, however, regrettably, the method proved ineffective in generating any words for $n\geq 3$, even when considering the partial intersection $\cap_{i=1}^n R_i$.

\begin{remark}
    We have also considered posibility of hyperbolic embedding, i.e. embeddings of the (compact region of) Cayley graph of $F$ into hyperbolic plane $\mathbb H$. We considered embedding by Sarkar algorithm~\cite{sarkar2012low} or as an orbit of $i$ under $\rho: F \hookrightarrow SL_2(\mathbb R)$ acting on the upper half-plane $\mathbb C_{\mathrm{Im}\geq 0}$ by Mobius transformations.

Both embeddings, although not isometric, retain certain metric properties of $F$, where the metric structure of $F$ is determined by the \emph{word metric}. However, during our preliminary investigations, we discovered that the metric structure of $F$ might not be the most suitable setting for optimization in our specific problem.
\end{remark}
 
\paragraph{Activation maximization}
While experimenting with neural based classifiers for elements in $R_i$, we investigated the feasibility of methods from the area of interpretable neural networks. Technique that we were studying, described in~\cite{nguyen2019understanding}, involves maximizing the activation of specific neurons in a neural network to generate images that are most likely to activate those neurons. We tried to adopt this method to the free group's word generation as it was done for protein sequences 
in~\cite{linder2021fast}.

Following~\cite{killoran2017generating}, \cite{linder2021fast}, we implemented the following approach. If the classifier input is given by sequence of length $M$ over the alphabet of size $N$, then the latent variable $z$ of size $N \times M$ is representing the probabilities of each token on each position. Suppose that classifier is trained up to a sufficiently high accuracy, then sampling with the probabilities from $z$ and feeding back the resulting words to the classifier, the relaxation procedure can be carried out in the latent space of $z$.

The appeal of this approach comes from the fact that it can be extended to an ensemble of classifiers. By maximizing the activation of all classifiers within the ensemble simultaneously using the same latent variable $z$, the approach can be used to generate words from the intersection $\cap_i R_i$ that maximize the activation of output neurons of all classifiers in the ensemble.

The appeal of this approach lies in its potential for extension to an ensemble of classifiers. By maximizing the activation of all classifiers within the ensemble simultaneously using the same latent variable $z$, the approach can be effectively applied across the entire ensemble. 

In the end, the activation maximization method turned out to be incapable of generating at even short non-trivial words from a single $R_i$. Despite the continuous effort, we were unable to generate words from the intersection using the joint activation maximization of multiple generators.

\paragraph{SeqGAN}
Another method that has been explored is generative adversarial networks (GANs)~\cite{jabbar2020survey}, specifically the SeqGAN~\cite{yu2017seqgan} approach. GANs are a type of deep neural network that can generate realistic synthetic data, such as images or text. SeqGAN focuses on generating realistic sequences of data, such as sentences or musical compositions. However, while GANs have shown promising results in some applications~\cite{jabbar2020survey} and can be used in similar circumstances~\cite{zehui2019intersectgan}, they did not yield promising results in our problem.
\paragraph{LSTM Ensemble}
Before developing the general masking approach, based on multi-labels, as mentioned in Remark~\ref{remark: masking ensemble}, we investigated its particular realization based on the ensemble of deep models. Our initial experiments were based on an ensemble of LSTM cells~\cite{hochreiter1997long}, in which each generator $G_j$ is trained on the subgroup $R_j$ or some symmetric commutator subgroup.

In generating the word from the intersection, to sample the next token from the ensemble we used a \emph{temperature sampling}, i.e. the probability of the next token $y_k$ from a sequence $y$ have a distribution

  \[
  y_k \sim \mathrm{Softmax}\left(\frac{1}{\tau}\sum\limits_{i = 1}^{N} w_i f_{\theta_i}(y)\right),
  \]
with weights $w_i$ for each generator $f_{\theta_i}$.

This approach allowed us to get first non-trivial elements from the full intersection. However, it was subsequently surpassed by a more general masking approach, which has also been demonstrated to be more scalable.

\subsection{Properties of sampling from $R_i$}
\begin{figure}[!htb]
\begin{subfigure}[b]{\textwidth}
    \includegraphics[width = 0.3\textwidth]{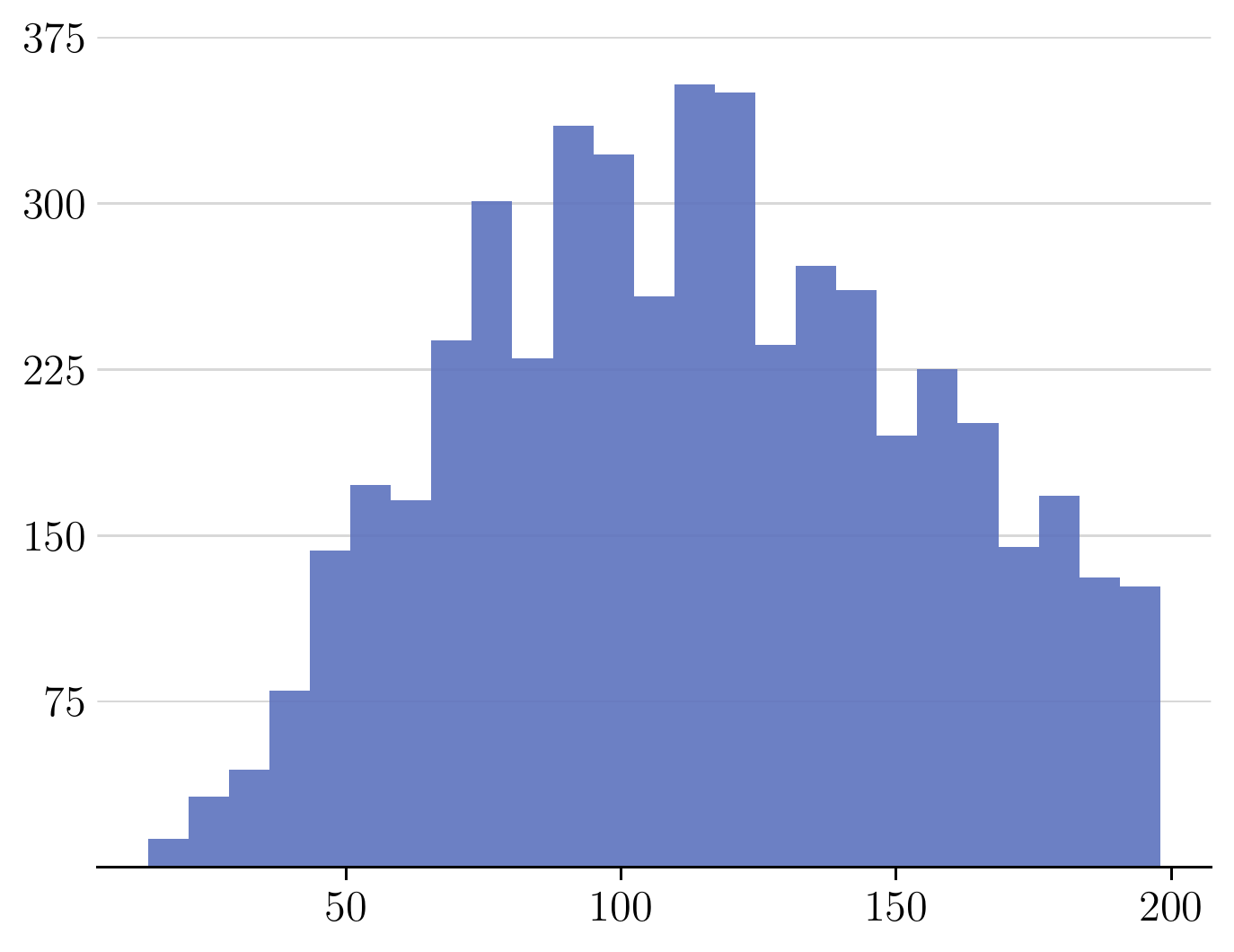}
    \includegraphics[width = 0.3\textwidth]{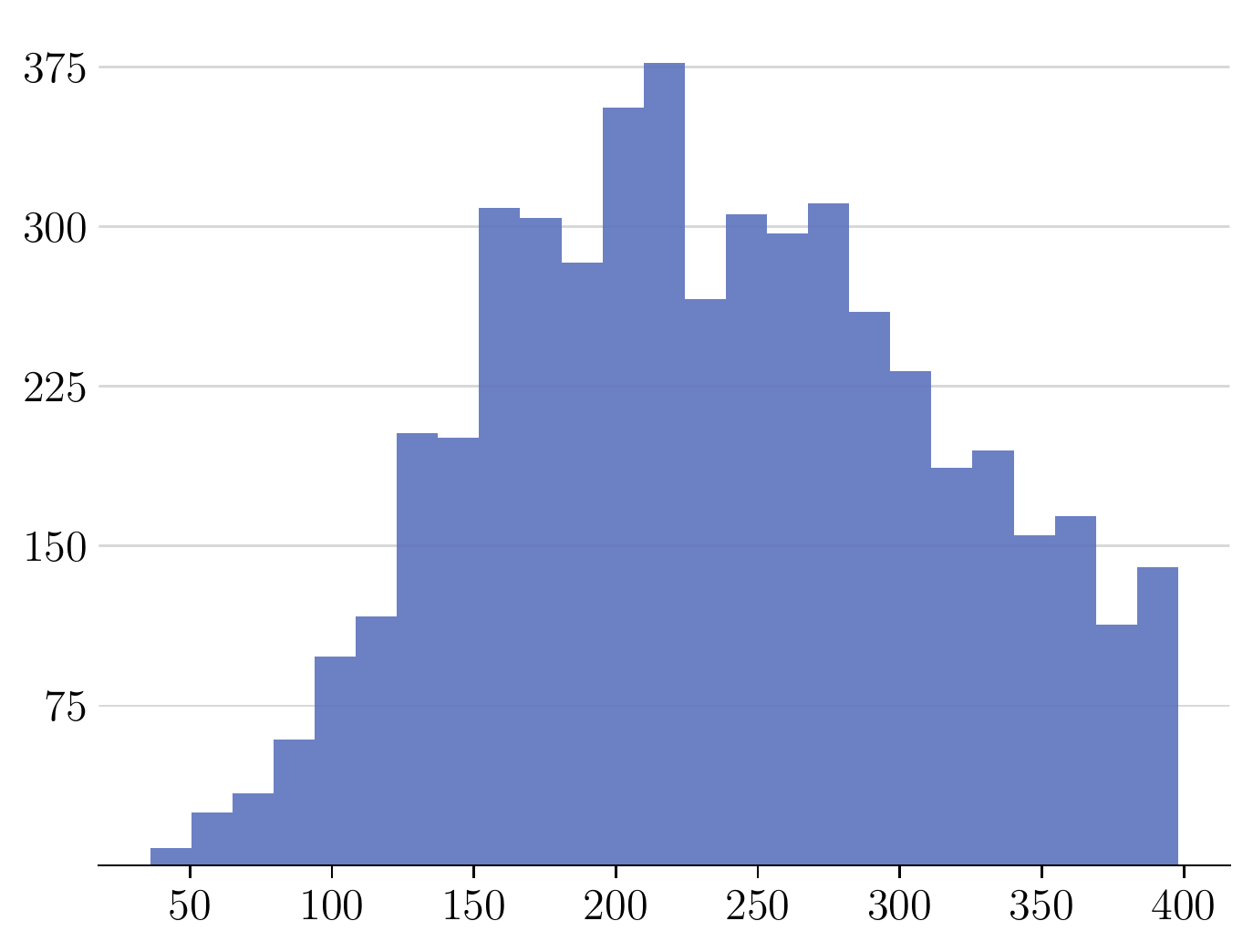}
    \includegraphics[width = 0.3\textwidth]{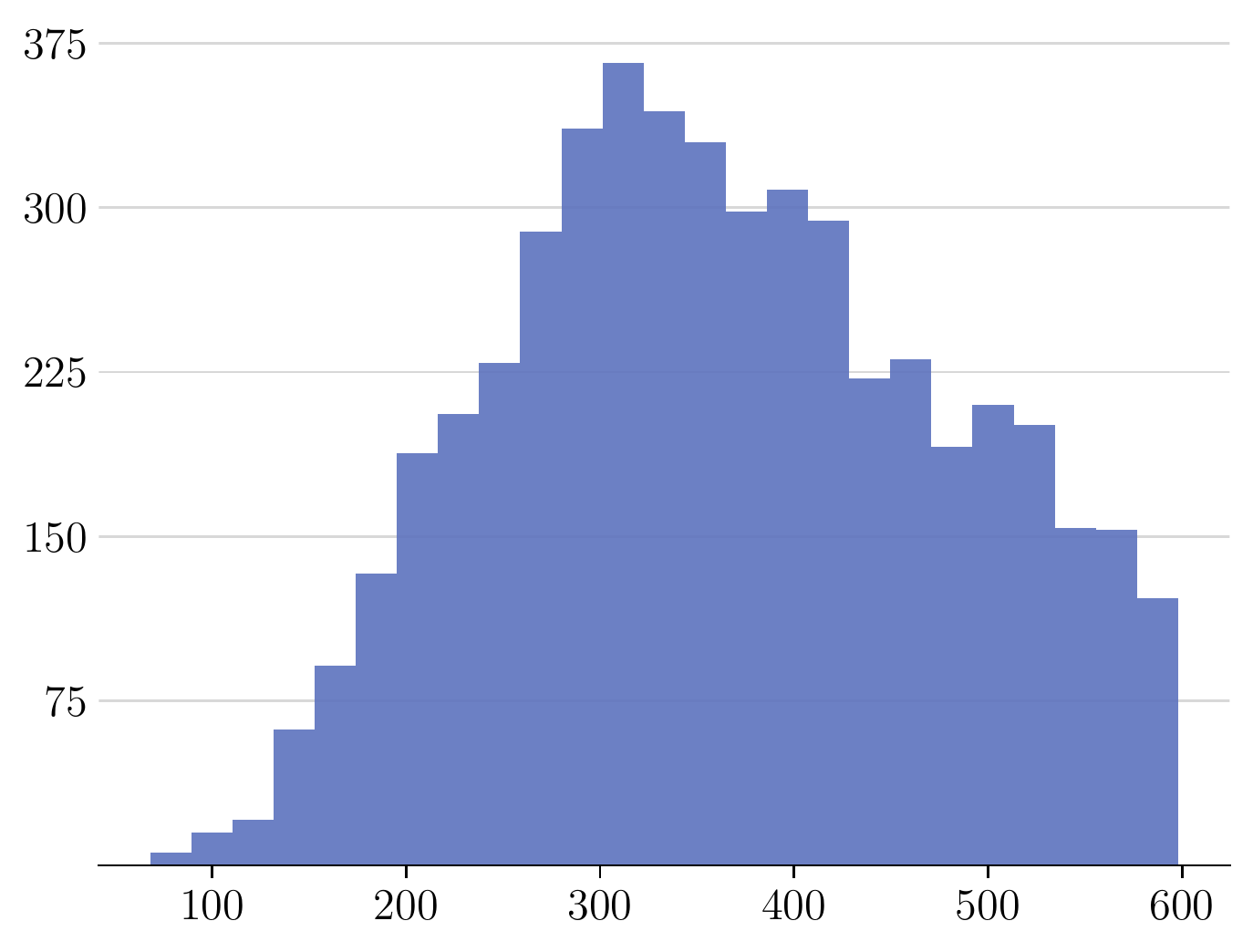}
\end{subfigure}
\begin{subfigure}[b]{\textwidth}
    \includegraphics[width = 0.3\textwidth]{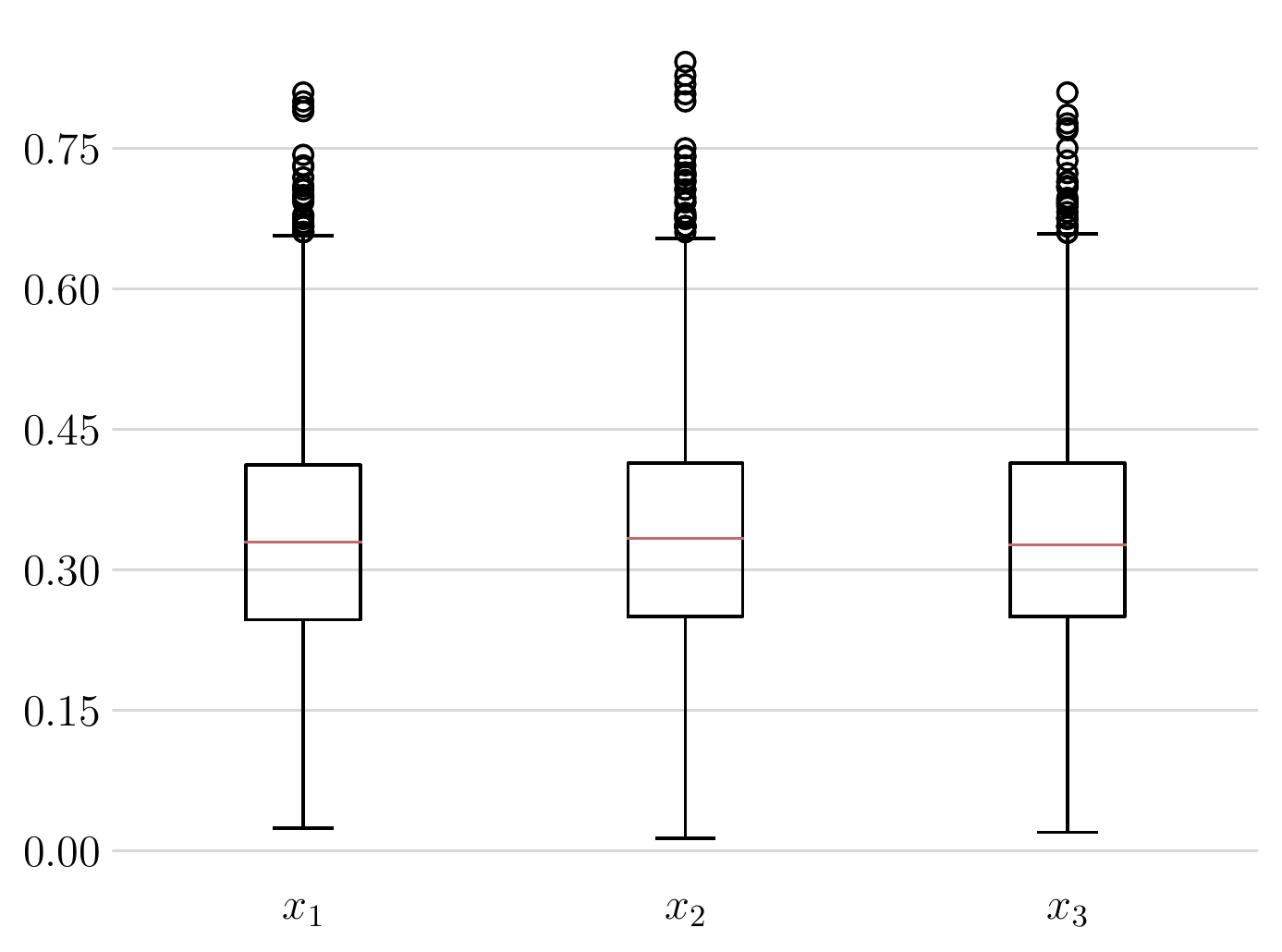}
    \includegraphics[width = 0.3\textwidth]{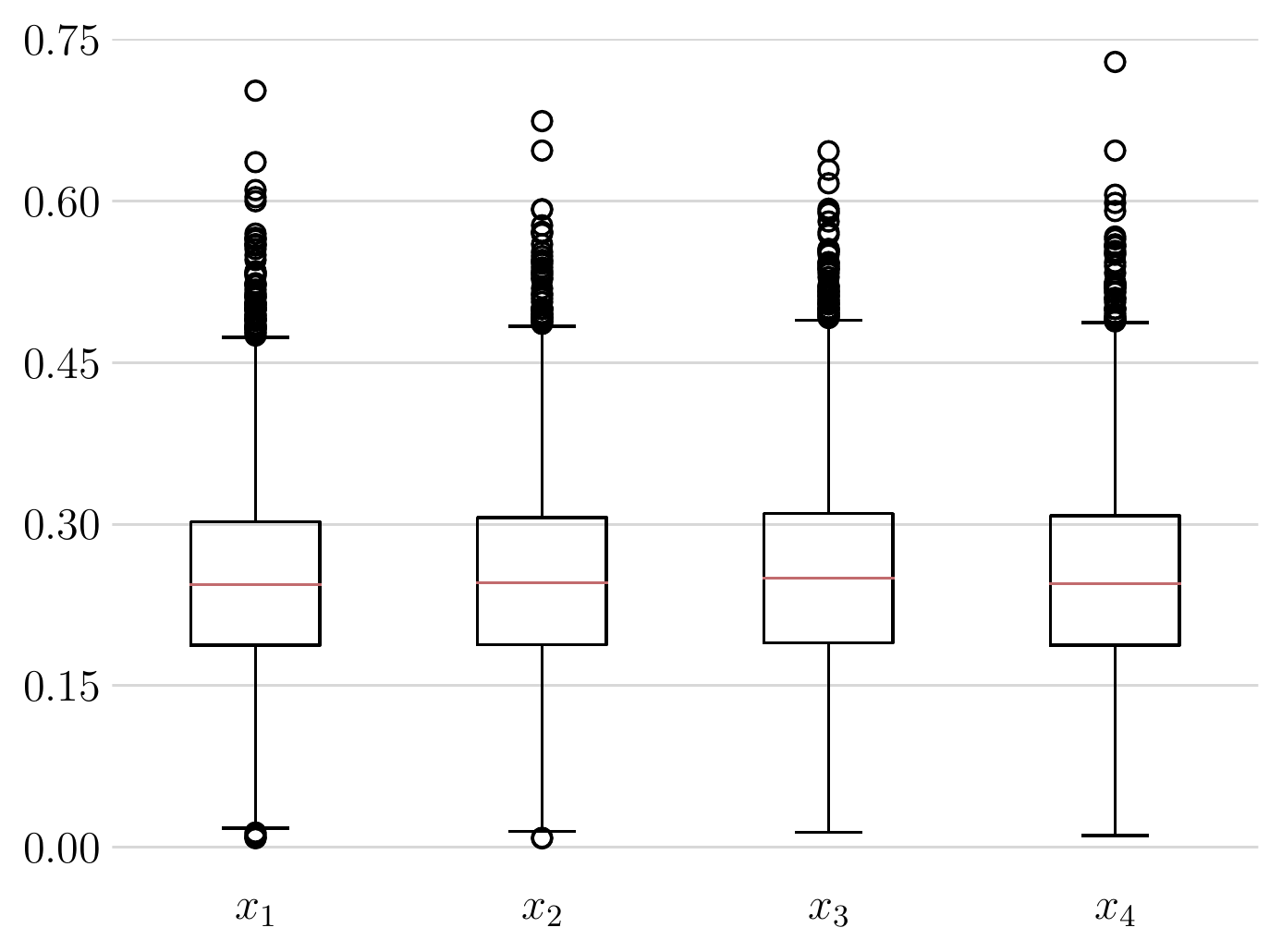}
    \includegraphics[width = 0.3\textwidth]{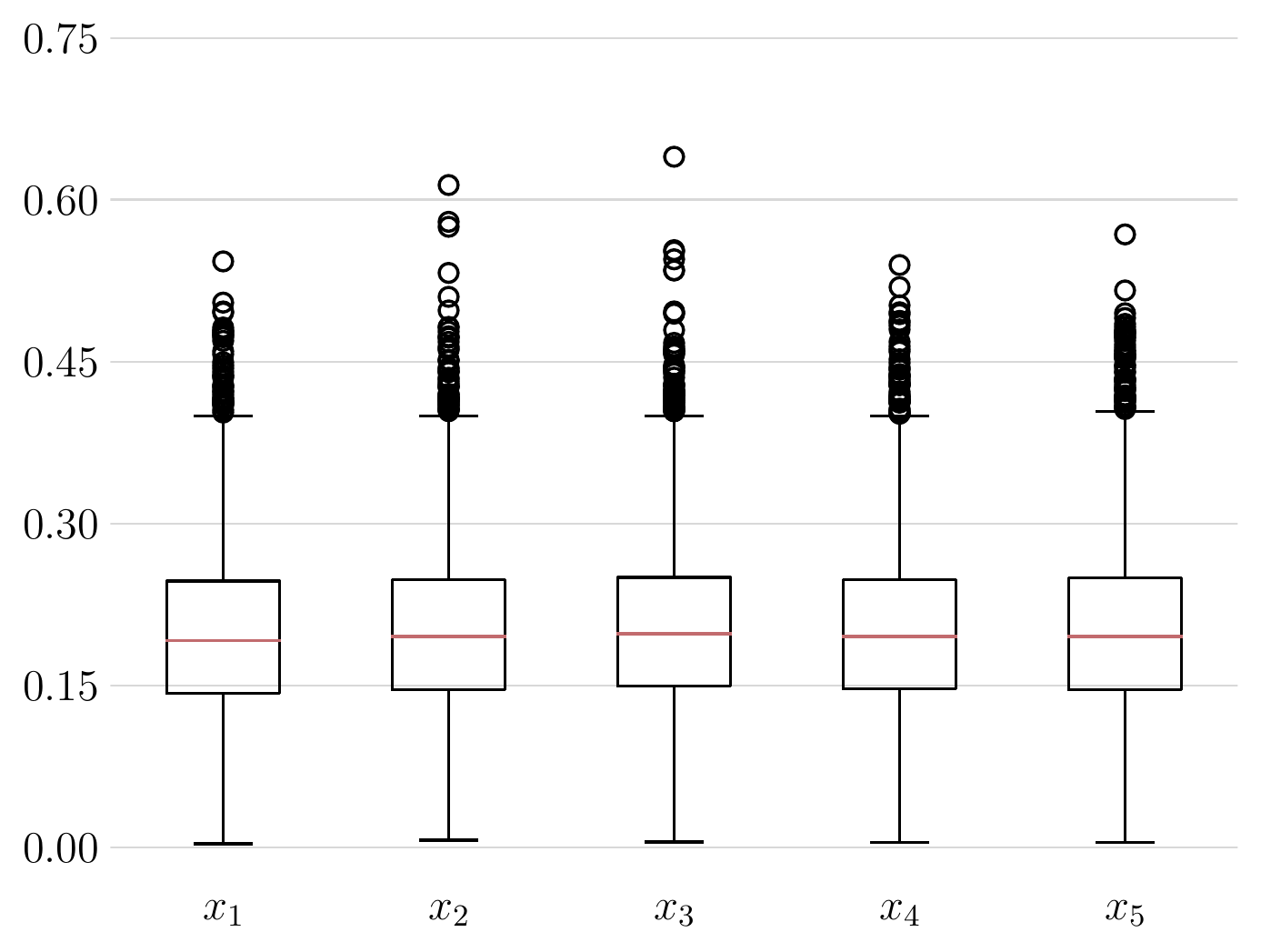}
\end{subfigure}

\caption{Length distribution (top) and occurrence ratio of generators $x_i$ (bottom) in the training dataset subsample for the $n = 3,4,5$. }
\label{fig:3-dataset-metrics}
\end{figure}

\subsection{Translator model}

For a given word $w \in [F,F]$ there are multiple way to present it in form of a product of commutators. Although there are algorithms to find such presentations~\cite{hall1950basis}, \cite{bartholdi2015commutator}, the presentations of long words in terms of basic commutators or as a product of ordinary (not iterated) commutators are bulk and practically unusable. We are looking for short and concise presentations, such as given in Table~\ref{table: non-trivial elements}.

To tackle this challenge, we employed an encoder-decoder model~\cite{vaswani2017attention} as a translator from a reduced word presentation to the commutator presentation.

The following procedure was implemented to generate the training and validation datasets:
\begin{enumerate}
\item A collection of random binary trees with a random number of leaves is generated.
\item The leaves of these trees are replaced with random words from a free group, ensuring that the word length remained within predefined bounds.
\item These modified trees (represented as strings) served as the labels for our model , representing a random commutators that we aim to translate.
\item Inputs, corresponding to the labels are given by reduced words, obtained from labels by opening all commutator brackets.
\end{enumerate}

Subsequently, we employed a standard training method for the translation model. We used BERT~\cite{devlin2019bert} with $2 \cdot 10^6$ parameters as an encoder as well as a decoder. We used AdamW with the learning rate $8\cdot 10^{-5}$ and the batch size equal to $64$. 

\begin{figure}
\includegraphics[width=0.5\textwidth]{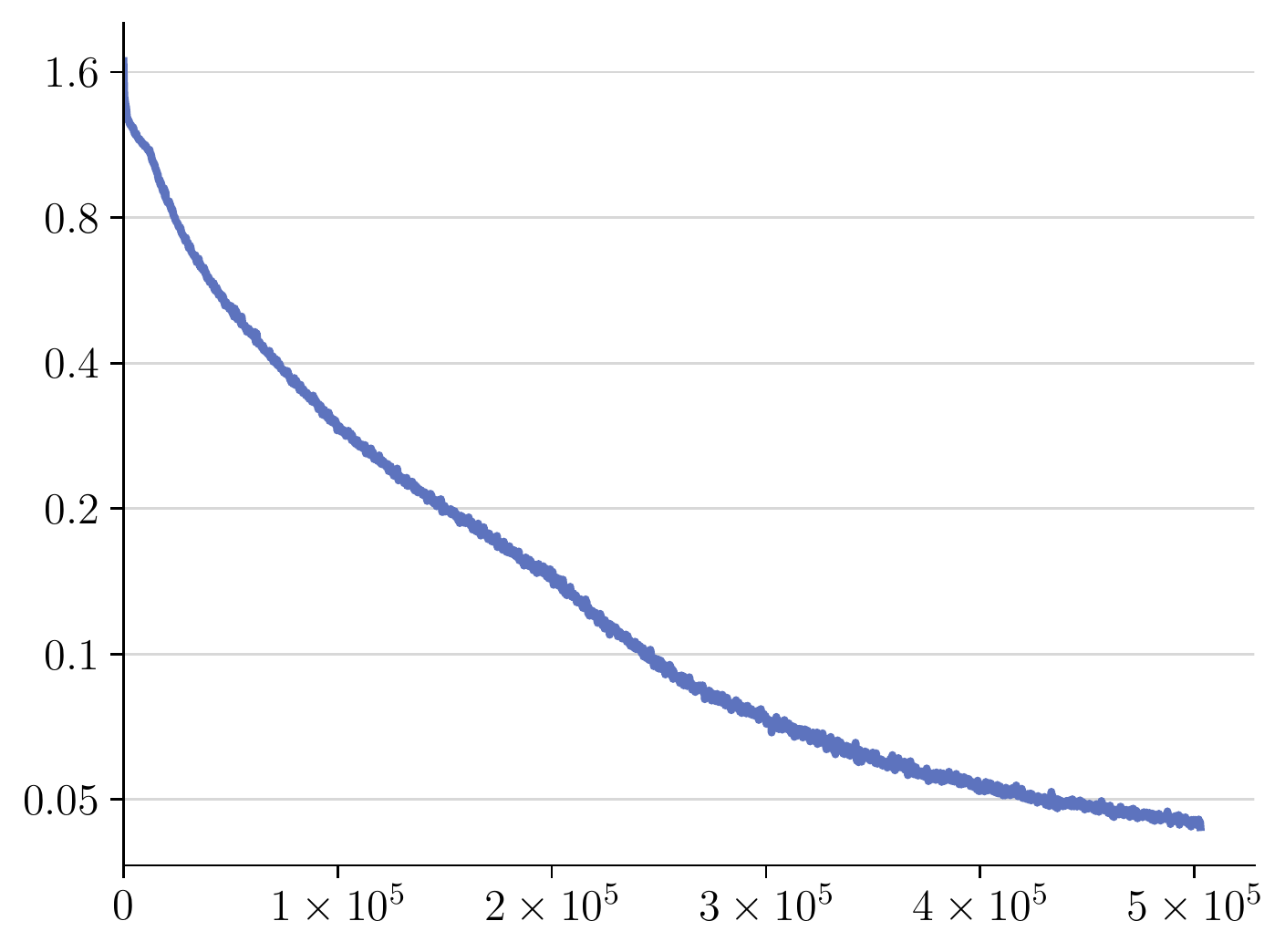}
\includegraphics[width=0.5\textwidth]{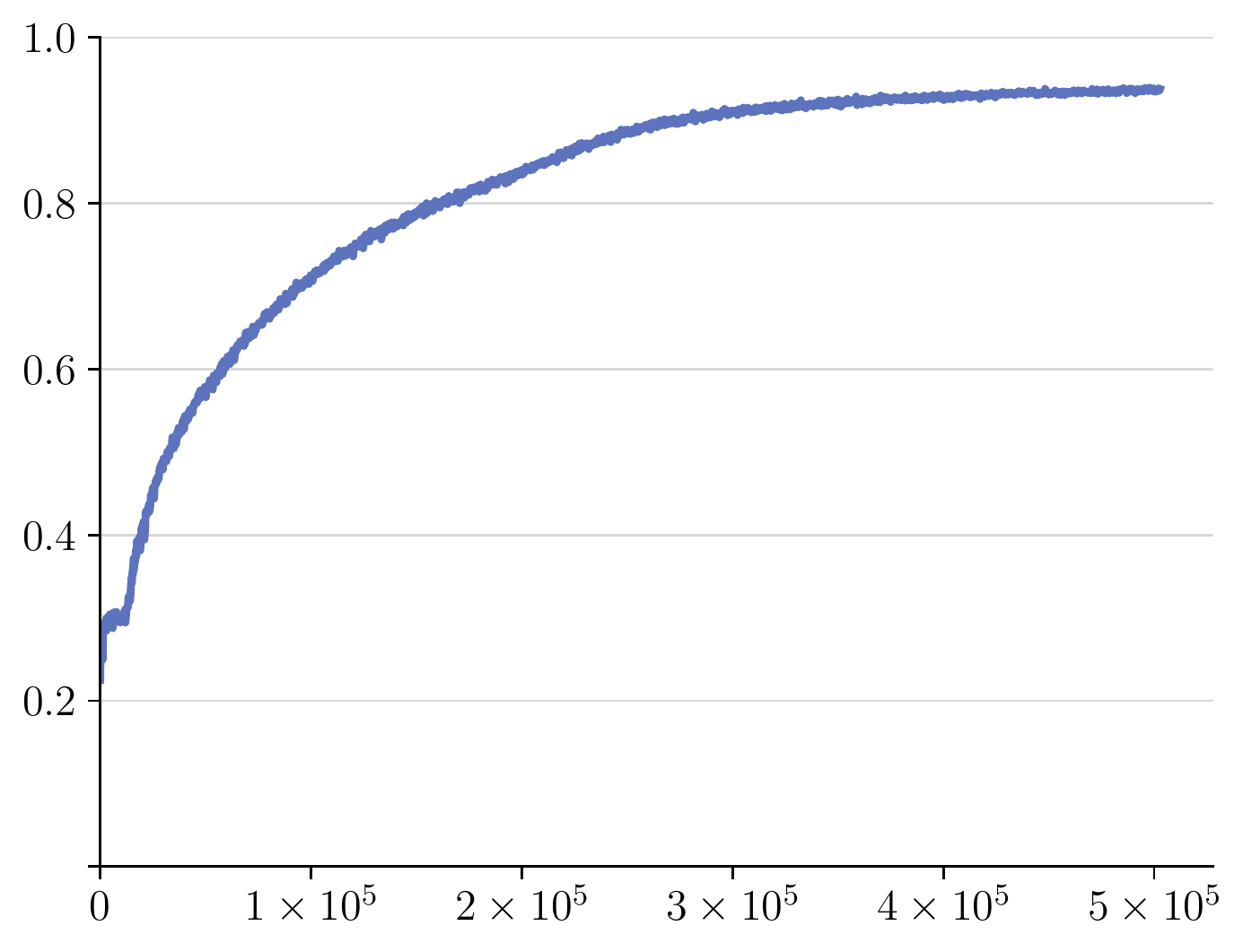}
\caption{Illustration of commutator translator learning process. Loss  with logarithmic scale (left) and BLeU score on validation (right).}
\label{fig: commutator translator learning}
\end{figure}

As with the main models (see Section~\ref{paragraph: datasets}), the training and validation datasets are generated in online mode with the following parameters:
\begin{itemize}
    \item number of generators $n = 4$,
    \item maximal length of a word = $250$,
    \item maximal length of a leave = $3$,
    \item maximal depth of a tree = $7$.
\end{itemize}
The training process is presented on Figure~\ref{fig: commutator translator learning}

\subsection{Examples of generated words}
We provide some examples of generated words for $n = 3,4$ in commutator presentation for various deep models and inference methods, see Section~\ref{subsec: model description} for additional details. The presentation is obtained using translator model described above.

\begin{itemize}

\item $n = 3$
\begin{gather*}
[[y^2,[y^{-1}x^{-1}z^{-1},[z^2,x]]],x]    \tag{prompt + beam}\\
\begin{rcases}
[[x,y],[yz,x]] \\
[[xy^{-2}x^{-1}z,y],[z,xy]] \\
[[[z,xy],y^{-1}xyx],y] 
\end{rcases}\tag{ignore + sample} \\
[[x, z], [x^{-1}, y][y, z]] \tag{mask + sample}
\end{gather*}

Note that some of the generated words are not in the symmetric commutator subgroup.

\item $n = 4$
\begin{gather*}
\begin{rcases}
[[xyzp,[y^4,[x,z]]],p]\\
[[(xyzp)^3,[[p^{-1},z^2],y]],x]  
\end{rcases}\tag{ignore + beam}\\
\begin{rcases}
[[yzpx^{-3},x^{-1}],[z^2,[p^2ypy^{-1},y^4]]]\\
[[[y,p^{-1}zp],xzp^2z^{-1}x^{-1}p],[x^3,zyzpxpxyx^{-1}]]   
\end{rcases}\tag{ignore + sample}\\
\begin{rcases}
[[[[z^{-1}y^{-1}x^{-1},p],y],x],z^3]\\
[[z^3,[(xyzp)^2,[y^2,x]]],p^3]   
\end{rcases}\tag{mask + beam}\\
[[p^2x^{-1}pxpx^{-1}px,[pzp^{-1}yzy^{-1},[y,x]]],xyzp]\tag{mask + sample}\\
\begin{rcases}
[[[[p^{-1}z^{-1}y^{-1}x^{-1},z],y],x],p]\\
[[x^3,[[x^{-1},p^{-1}z^{-1}y^{-1}],[p,y]]],z^2]   
\end{rcases}\tag{prompt + beam}\\
\begin{rcases}
[[[p^{-1}y^2p,z],p^{-1}xz^{-1}xzp],[xyzp,zp^3xpx^{-1}z^{-1}]]\\
[[[zpx,y],x^{-1}zx^2zx^{-1}z^{-1}],[[x^3yx,y^{-1}],p^3]]   
\end{rcases}\tag{prompt + sample}\\
\end{gather*}

\end{itemize}

\subsection{Details on greedy algorithm}
Suppose we have a given prefix $p = xyz$ and we aim to generate a word from the intersection of the normal closures $R_1 = \langle x \rangle^F$, $R_2 = \langle y \rangle^F$, and $R_3 = \langle z \rangle^F$. The corresponding stacks $S_i$ created for each normal closure $R_i$ are $S_1 = [y,z]$, $S_2 = [x,z]$, $S_3 = [x,y]$. In a list of candidates for the next token, the token $z^{-1}$ is excluded because it would reduce the length of the prefix $p$.

Next, the algorithm assigns points to each token based on the two criteria mentioned in Section~\ref{paragraph: greedy}. The token $y^{-1}$ would receive a point because it reduces the length of the $S_3$. The tokens $x$ and $x^{-1}$ would receive points because they do not increase the length of $S_1$. Similarly, $y$, $y^{-1}$, and $z$ would receive points because they do not increase the length of $S_2$ and $S_3$ respectively. The algorithm then selects the token with the highest point value as the next token to add to the prefix. In this case, $y^{-1}$ has the highest number of points, so it would be selected as the next token to add to the prefix.

In the case of a ``many-token'' normal closure, such as $R_0 = \langle xyz \rangle^F$, tokens receive points based on their ability to bring the top of the corresponding stack $S_0$ closer to a rotation of the normal closure. In particular, if token $x_k$ is the top of $S_0$, then the token $x_{(k+1) \textrm{ mod } (n+1)}$ will receive a point.
\newpage
\subsection{Auxilary learning curves}
\begin{figure*}[!htb]
     \centering
     \begin{subfigure}[b]{0.32\textwidth}
         \centering
         \includegraphics[trim={0.25cm 0.3cm 0.25cm 0.3cm},clip,width = \textwidth]{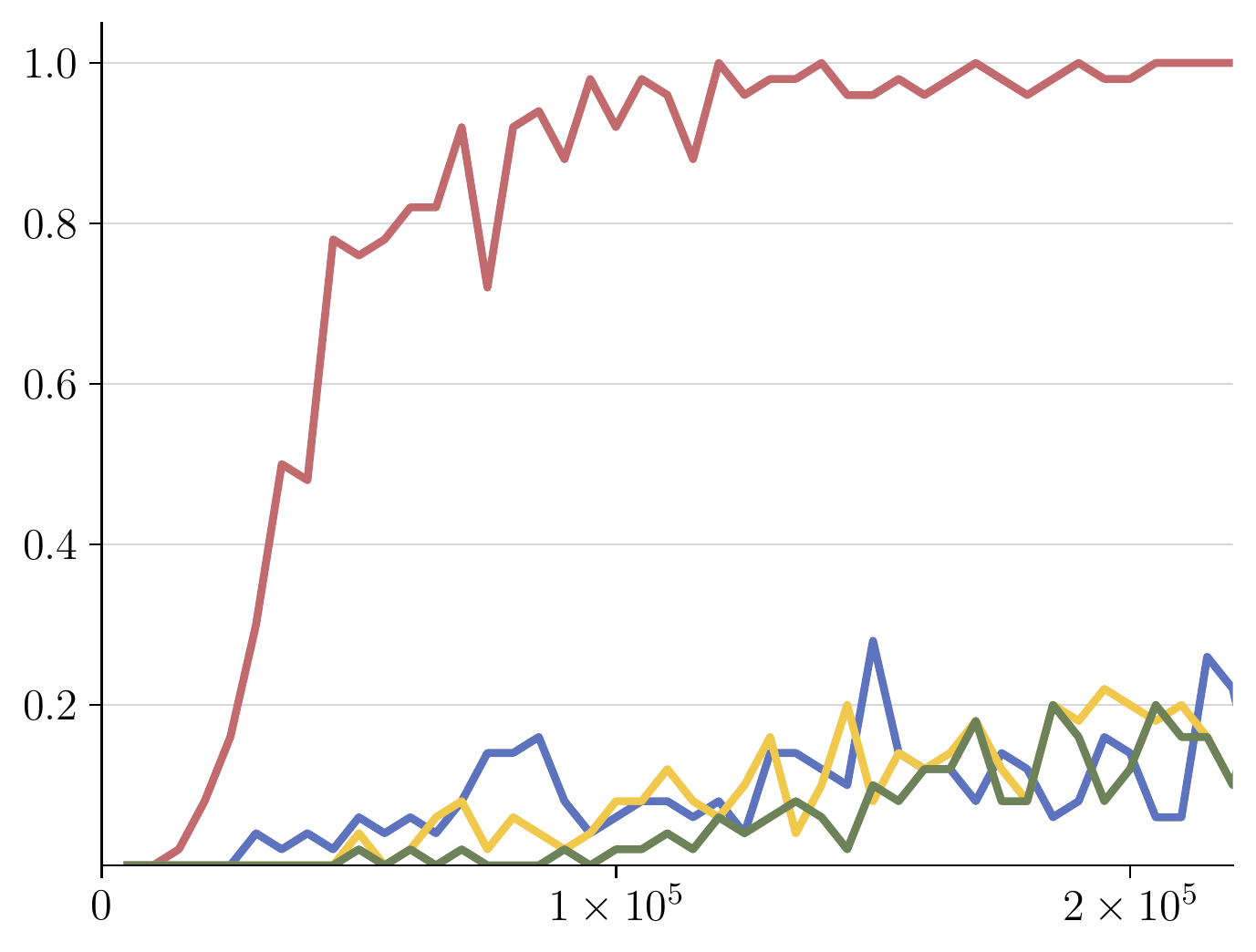}
         \caption{$k = 5$}
     \end{subfigure}
     \hfill
     \begin{subfigure}[b]{0.32\textwidth}
         \centering
         \includegraphics[trim={0.25cm 0.3cm 0.25cm 0.3cm},clip,width=\textwidth]{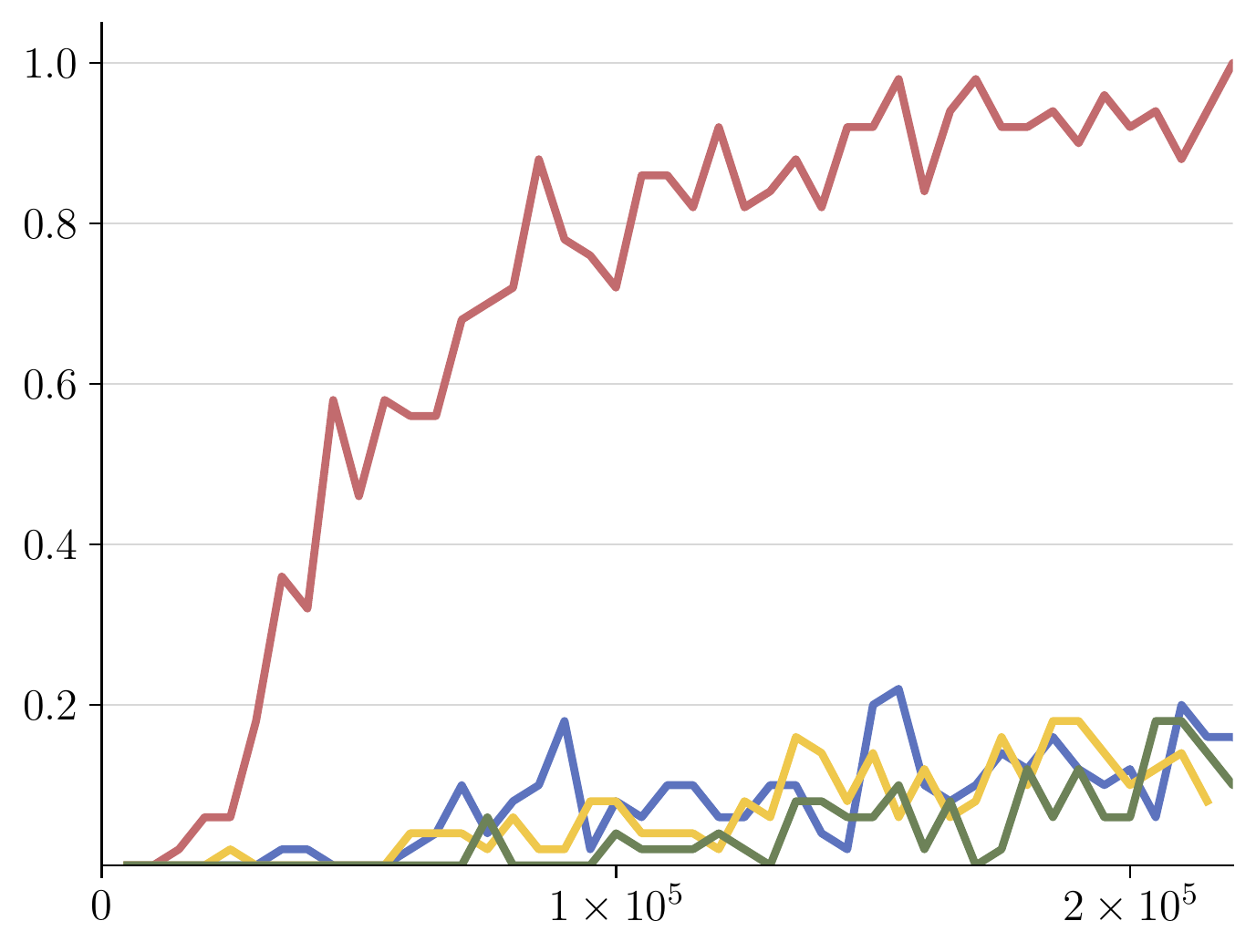}
         \caption{$k = 7$}
     \end{subfigure}
     \hfill
     \begin{subfigure}[b]{0.32\textwidth}
         \centering
         \includegraphics[trim={0.25cm 0.3cm 0.25cm 0.3cm},clip,width=\textwidth]{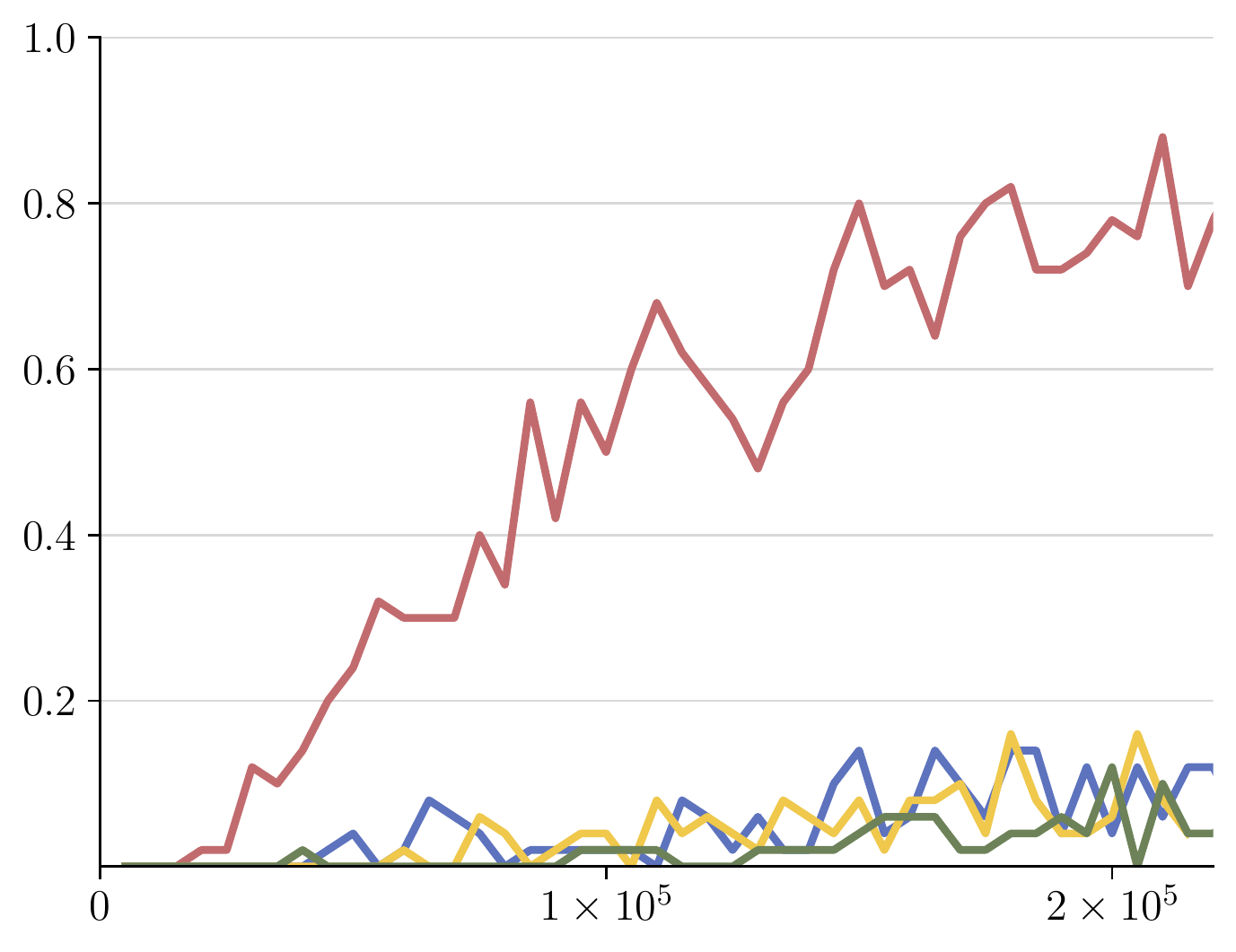}
         \caption{$k = 10$}
     \end{subfigure}
        \caption{Completion ratio  with nucleus sampling for $n=3$ and various length of prefixes $k$.}
        \label{fig: aux learning curve n = 3 sample}
\end{figure*}

\begin{figure*}[!htb]
     \centering
     \begin{subfigure}[b]{0.32\textwidth}
         \centering
         \includegraphics[trim={0.25cm 0.3cm 0.25cm 0.3cm},clip,width = \textwidth]{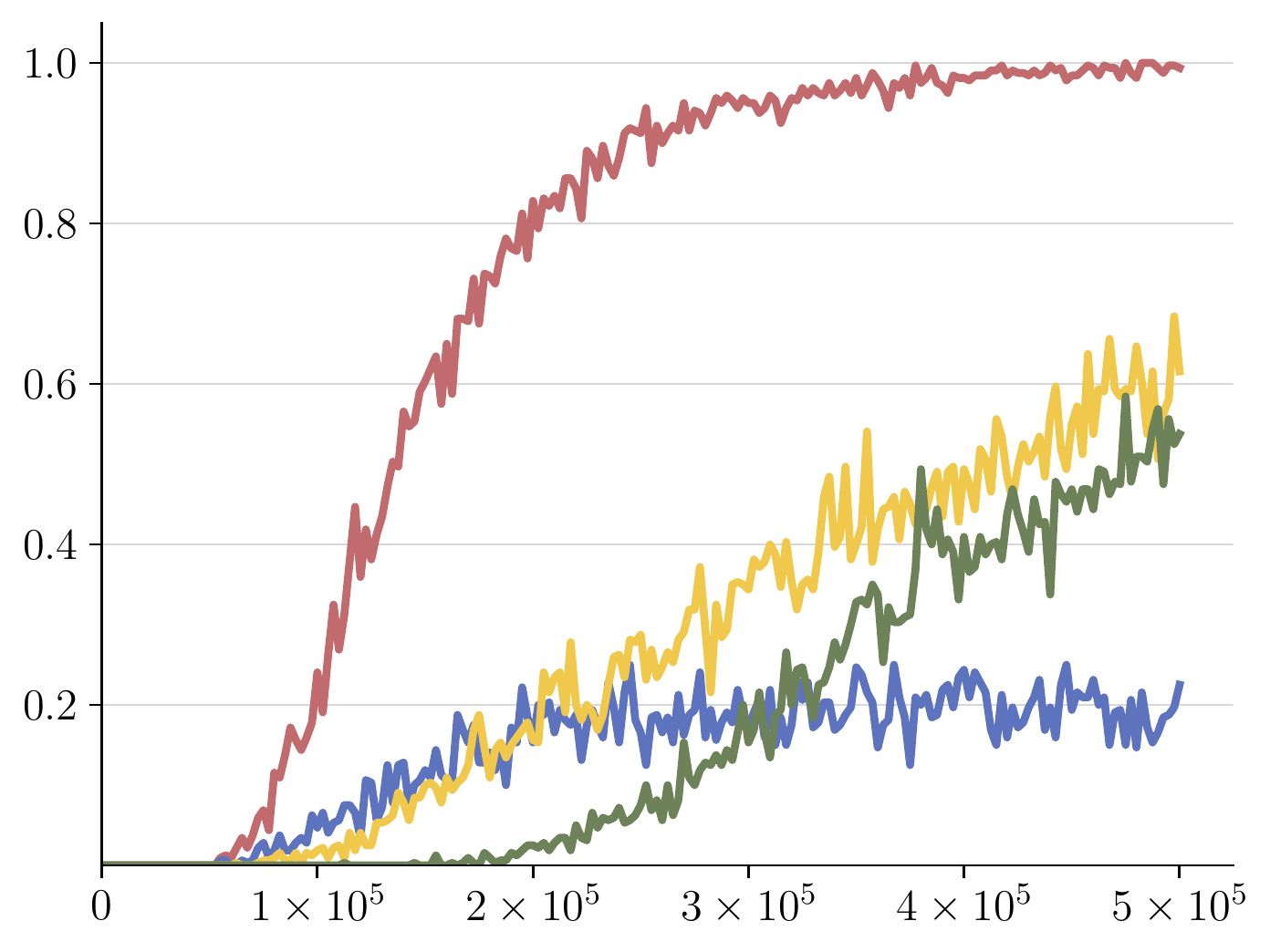}
         \caption{$k = 5$}
     \end{subfigure}
     \hfill
     \begin{subfigure}[b]{0.32\textwidth}
         \centering
         \includegraphics[trim={0.25cm 0.3cm 0.25cm 0.3cm},clip,width=\textwidth]{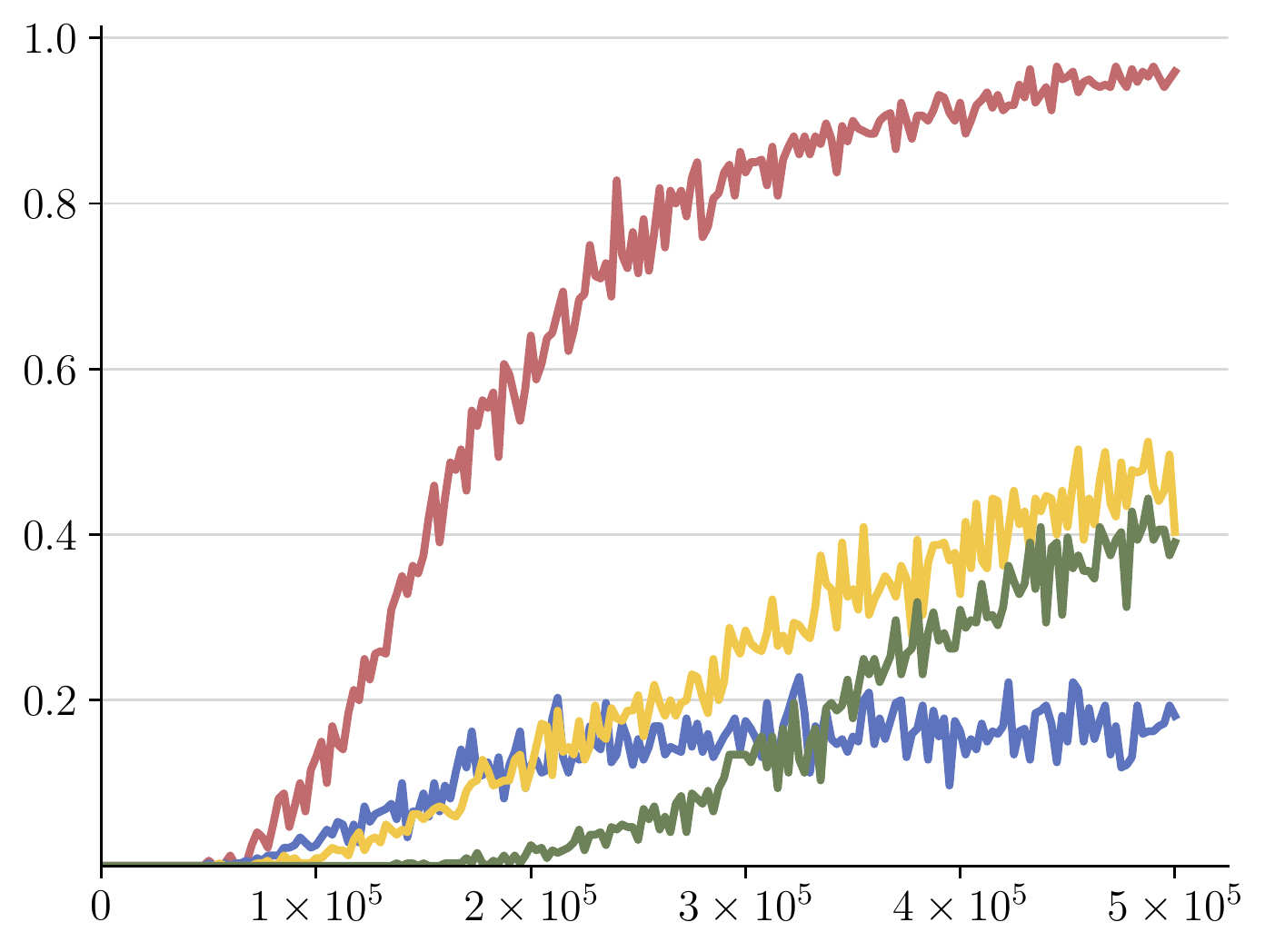}
         \caption{$k = 7$}
     \end{subfigure}
     \hfill
     \begin{subfigure}[b]{0.32\textwidth}
         \centering
         \includegraphics[trim={0.25cm 0.3cm 0.25cm 0.3cm},clip,width=\textwidth]{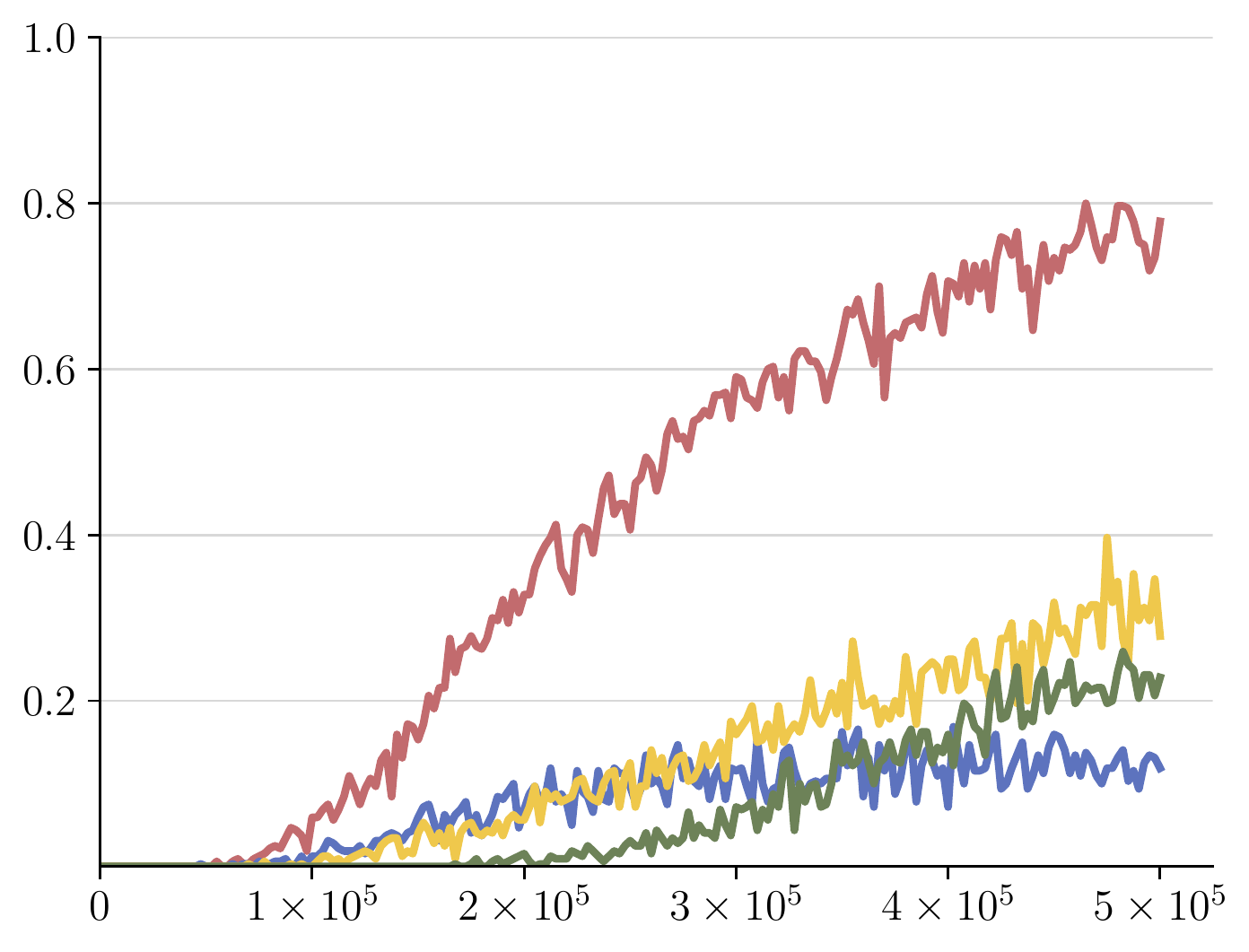}
         \caption{$k = 10$}
     \end{subfigure}
        \caption{Completion ratio  with nucleus sampling for $n=4$ and various length of prefixes $k$.}
        \label{fig: aux learning curve n = 4 sample}
\end{figure*}

\begin{figure*}[!htb]
     \centering
     \begin{subfigure}[b]{0.32\textwidth}
         \centering
         \includegraphics[trim={0.25cm 0.3cm 0.25cm 0.3cm},clip,width = \textwidth]{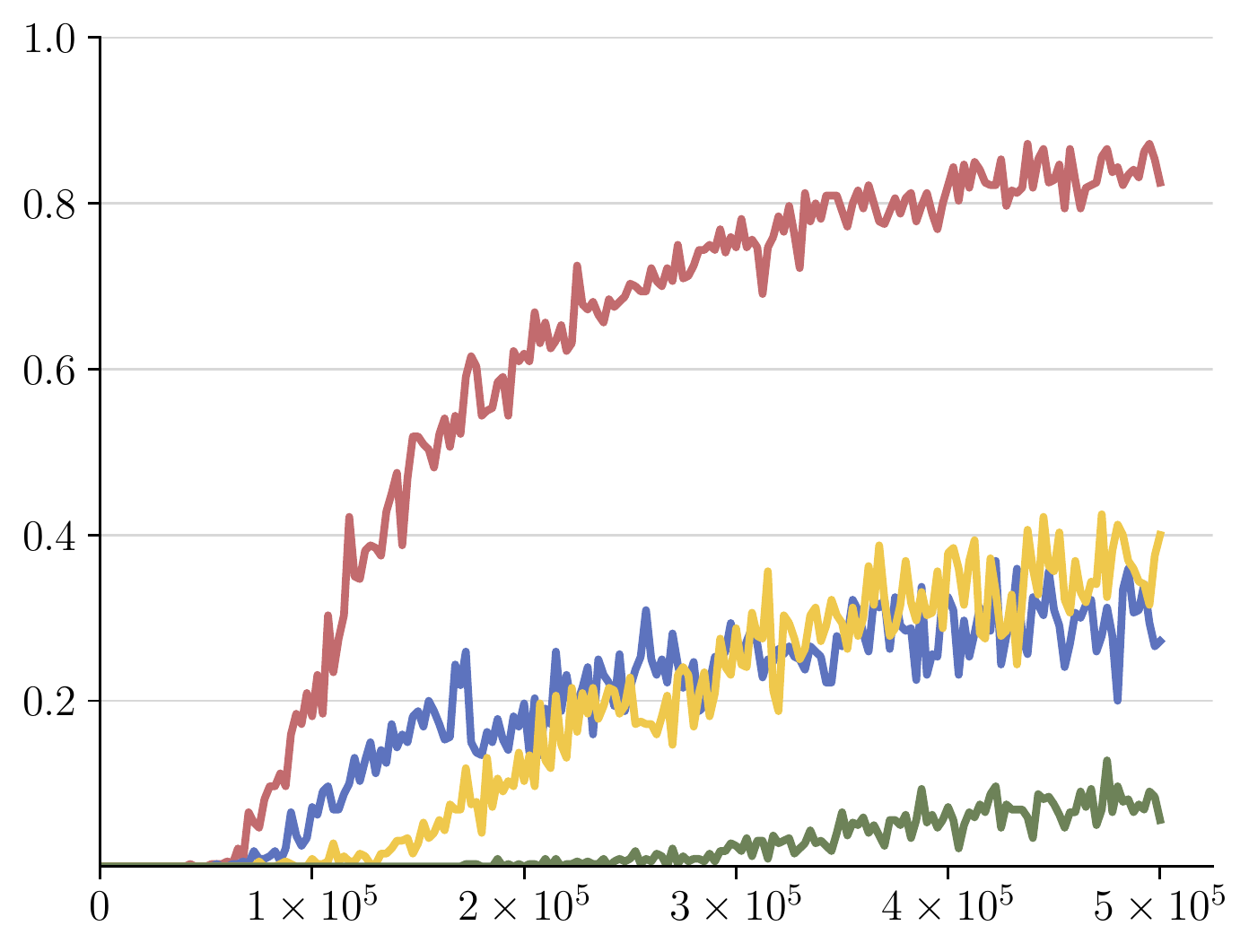}
         \caption{$k = 5$}
     \end{subfigure}
     \hfill
     \begin{subfigure}[b]{0.32\textwidth}
         \centering
         \includegraphics[trim={0.25cm 0.3cm 0.25cm 0.3cm},clip,width=\textwidth]{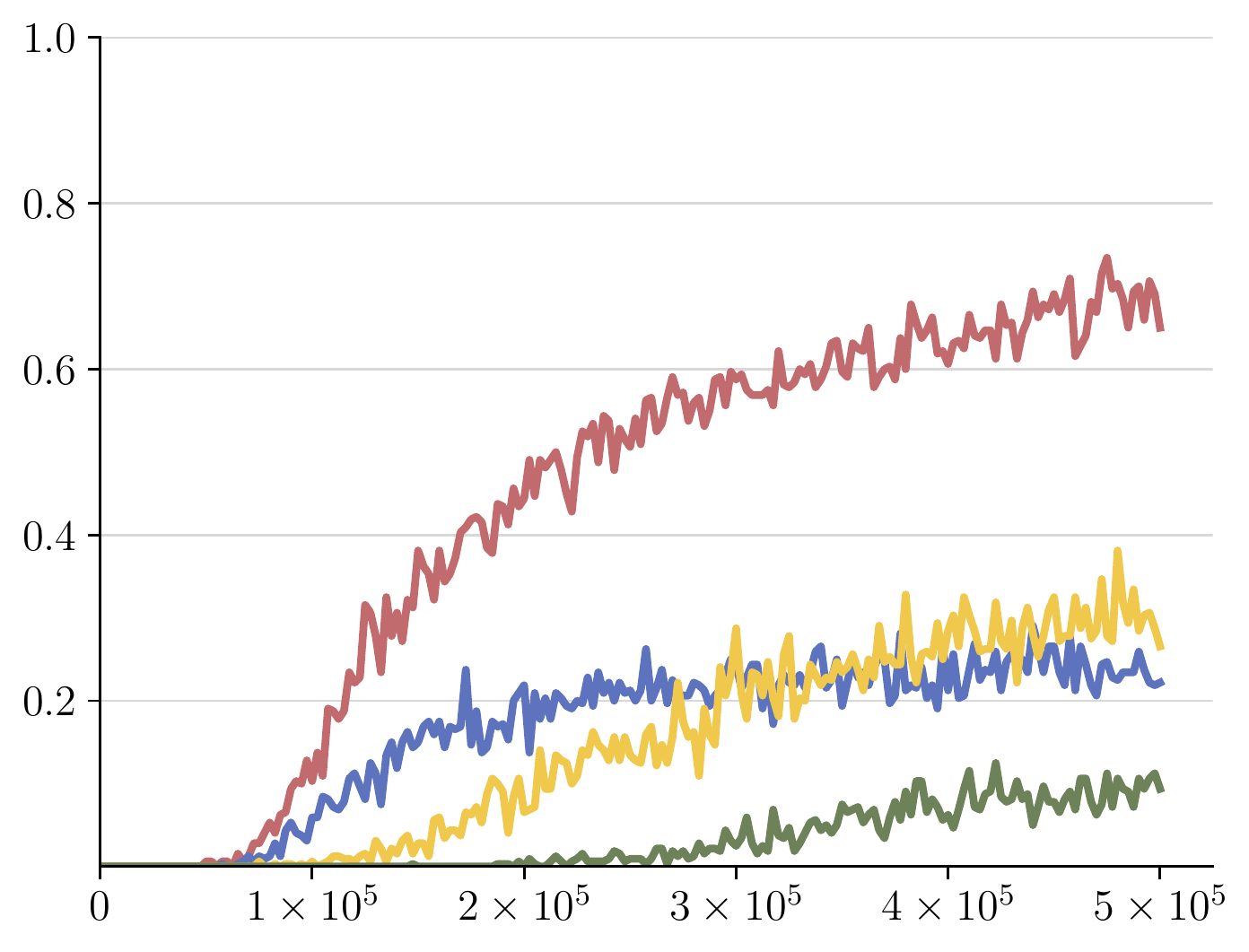}
         \caption{$k = 7$}
     \end{subfigure}
     \hfill
     \begin{subfigure}[b]{0.32\textwidth}
         \centering
         \includegraphics[trim={0.25cm 0.3cm 0.25cm 0.3cm},clip,width=\textwidth]{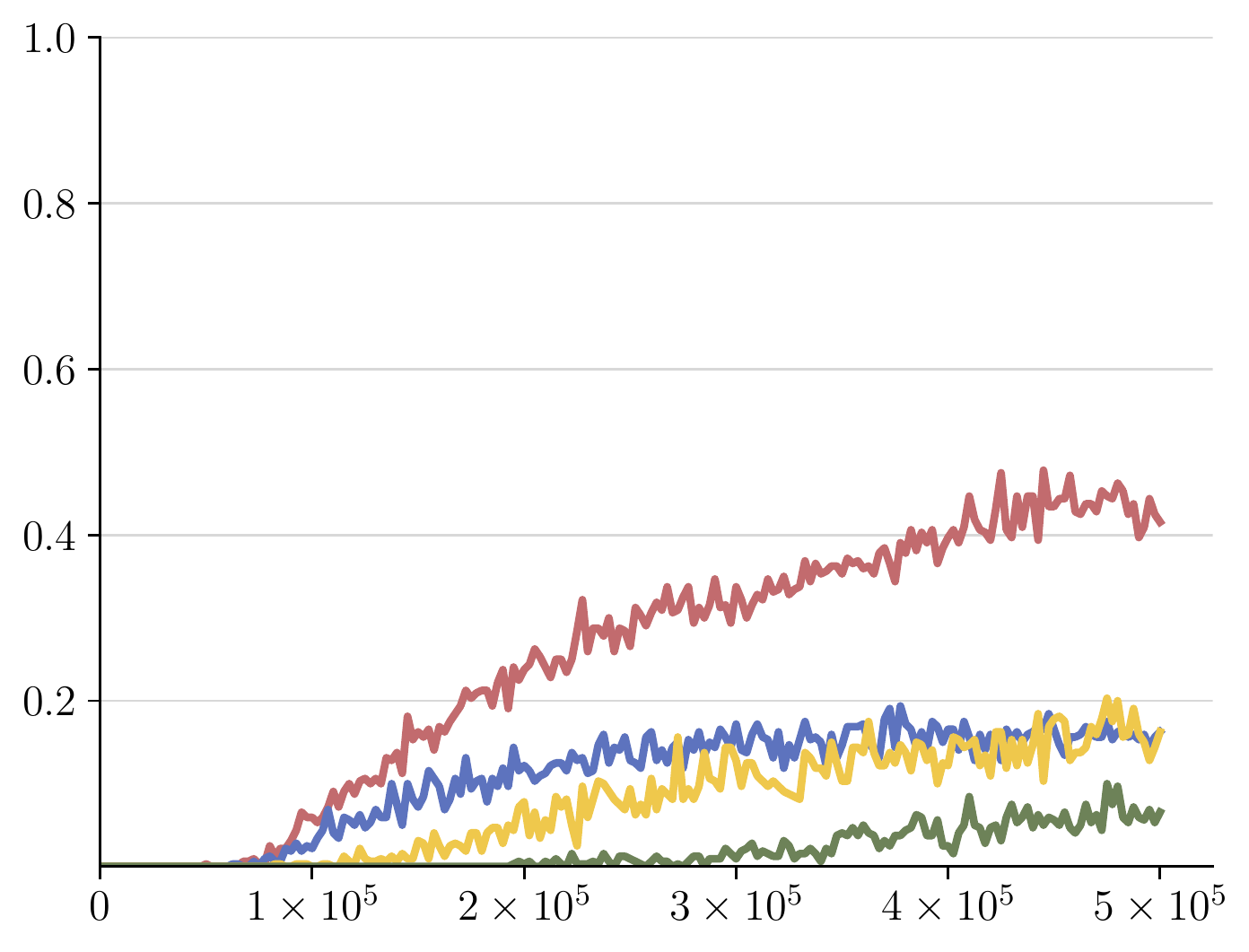}
         \caption{$k = 10$}
     \end{subfigure}
        \caption{Completion ratio  with beam search for $n=4$ and various length of prefixes $k$.}
        \label{fig: aux learning curve n = 4 beam}
\end{figure*}

In Figures~\ref{fig: aux learning curve n = 3 sample}-\ref{fig: aux learning curve n = 4 beam} we present additional graphs of learning process for various hyperparameters like inference method and prefix length. Note that variable prefix allows more variability in models output as well. Color coding is the same as in Figure~\ref{fig: Metrics}: red for \emph{negative baseline}, yellow for \emph{prompt}, green for \emph{masking}, blue for \emph{ignore}.

\end{document}